\documentclass[9pt,journal,twocolumn]{IEEEtran}

\IEEEoverridecommandlockouts                              
\usepackage{color}
\usepackage[dvipsnames]{xcolor}
\colorlet{LightRubineRed}{RubineRed!70!}
\colorlet{Mycolor1}{green!10!orange!90!}
\definecolor{Mycolor2}{HTML}{00F9DE}

\usepackage{calc}

\usepackage{tikz}
\tikzset{
sometest/.style={
    align=left,
    text width=5cm
},
asdfasdf/.style={
    rectangle
},
}
 \usepackage{subcaption}

\usepackage{longtable}
\usepackage{amsmath}
\usepackage{graphicx}
\usepackage{color}

\usepackage[dvips]{epsfig}
\usepackage{comment}
\usepackage{cite}
\usepackage{graphicx, amssymb}
\usepackage{amsmath}

\newcommand{\beq}{\begin{equation}}
\newcommand{\eeq}{\end{equation}}
\newcommand{\carre} {\hfill $\blacksquare$}

\newcommand{\startcompact}[1]{\par\vspace{-0.75em}\begin{#1}%
\allowdisplaybreaks\ignorespaces}

\newcommand{\stopcompact}[1]{\end{#1}\ignorespaces}

\tikzstyle{vertex}=[circle, draw,  thick, inner sep=0pt, minimum size=9pt,fontsize=9pt]

\newtheorem{theo}{Theorem}

\newtheorem{rem}{Remark}
\newtheorem{deff}{Definition}
\newtheorem{lem}{Lemma}
\newtheorem{corr}{Corollary}

\newtheorem{pro}{Proposition}

\title{\LARGE \bf Strong Structural Controllability of Networks under  Time-Invariant and Time-Varying Topological Perturbations \color{black}}

\author{Shima Sadat Mousavi, Mohammad Haeri, and Mehran Mesbahi
\thanks{S.S. Mousavi was with the Department of Electrical Engineering, Sharif University of Technology, Tehran, Iran. She is now with the Department of Civil, Environmental and Geomatic Engineering, ETH Z\"urich, 8093 Z\"urich, Switzerland (e-mail: shimasadat$\_$mousavi@ee.sharif.edu, mousavis@ethz.ch).}
\thanks{M. Haeri is with the Department of Electrical Engineering, Sharif University of Technology, Tehran, Iran (e-mail: haeri@sharif.ir).}%
\thanks{M. Mesbahi is with the Department of Aeronautics and Astronautics, University of Washington, WA 98195 (e-mail: mesbahi@uw.edu).}
\thanks{A preliminary version of this work  was presented at the 2017 American Control Conference, Seattle, May 24-26, 2017  \cite{mousavi2017robust}.}
}

\begin{document}

\maketitle
\thispagestyle{empty}
\pagestyle{empty}

\begin{abstract}
This paper investigates the robustness of strong structural controllability for linear time-invariant and linear time-varying directed networks with respect to structural perturbations, including edge deletions and additions. 
In this direction, we introduce a new construct referred to as a perfect graph associated with a network with a  given set of control nodes. The tight upper bounds on the number of edges that can be added to, or removed from a network, while ensuring strong structural controllability, are then derived. Moreover, we obtain a characterization of critical edge-sets, the maximal sets of edges whose any subset can be respectively added to, or removed from a network, while preserving strong structural controllability. In addition, procedures for combining networks to obtain  strongly structurally controllable network-of-networks are proposed. Finally, controllability conditions are proposed for networks whose edge weights, as well as their structures, can vary over time. 
\end{abstract}



\section{Introduction}

In recent years, controllability analysis of 
networks from a graph-theoretic point of view has become an active area of research 
~\cite{mousavi2018controllability,mousavi2018laplacian}. Among the various approaches adopted to reason about network controllability, notions of weak structural controllability (ws-controllability) and strong structural controllability (ss-controllability)  capture two facets of examining controllability for a parameterized family of linear time-invariant (LTI) systems. For both approaches, system parameters  
are classified into parameters that (always) assume a zero or nonzero value, and the exact value of nonzero parameters is 
unknown. In the weak structural framework, controllability results can be extended to {\em almost all}  networks with the same structure (see e.g., \cite{lin1974structural,liu2011controllability}), while in the strong structural setting that provides a stronger notion of controllability, the results hold for {\em all} networks with the same zero-nonzero pattern.
Mayeda and Yamada introduced the notion of ss-controllability for single-input LTI systems \cite{mayeda1979strong}. Their results were then extended to multi-input systems in \cite{bowden2012strong}. Subsequently, the notion of ss-controllability was further explored through such graph-theoretic concepts as cycle families and constrained matchings \cite{jarczyk2011strong,chapman2013strong}. In particular, \cite{monshizadeh2014zero} provided a necessary and sufficient condition for ss-controllability of an LTI network in terms of the notion of a zero forcing set (ZFS); \color{black}{this result was further extended in \cite{trefois2015zero,monshizadeh2015strong,mousavi2016controllability, van2017distance,mousavi2018structural,mousavi2018null}.  
 \color{black}Other references on ss-controllability of networks include \cite{reissig2014strong,gracy2017structural,jia2018sufficient,sadat2019strong}.}
\color{black} 

In this paper, using the notion of a ZFS for investigating ss-controllability,
 we provide new insights into the robustness of controllability to uncertainties in the network. It should be noted that the ss-controllability results are typically valid when no edges are added to, or removed from the network, i.e., the underlying pattern of zero/nonzero elements remains intact. \color{black}{The only structural changes, being allowed in some ss-controllability results, include  adding or removing self-loops \cite{monshizadeh2014zero,van2017distance}.} \color{black} However, since networks may experience \textit{structural perturbations} (e.g., cellular biochemical pathways in biological networks can be altered), robustness analysis for network controllability with respect to both uncertainties in the system parameters and presence or absence of interconnections becomes of paramount importance. Examples of such structural perturbations include loss of interconnections in a power distribution network caused by link failures (or due to malicious attacks) and failed or new social ties in a social network. Another case for robustness analysis with respect to structural perturbations is made when only an approximate representation of a network has been provided. For example, it is common to omit edges with small edge weights in a network. However, the existence of these edges can alter ss-controllability properties of the network. Recent works on robustness of  ws-controllability of networks against link failures include~\cite{rahimian2013structural, nie2014robustness, mengiste2015effect, lu2016attack, monnot2016sensitivity}.  We note that ws-controllability is not affected by edge additions, while removal of edges may affect it. On the other hand, both edge additions and deletions may change ss-controllability. \color{black}{ To the best of our knowledge, \cite{mousavi2017robust} is the only work in the literature, investigating the robustness of ss-controllability with respect to different types of structural perturbations; \color{black} the current paper extends the results of \cite{mousavi2017robust} to strongly structurally controllable (SSC) networks-of-networks and linear time-varying (LTV) networks.} \color{black}

\color{black}The notion of a ZFS, which is related to a particular coloring of nodes of a graph,  was first introduced in~\cite{work2008zero} to study the minimum rank problem for symmetric patterned matrices; the work \cite{barioli2009minimum} further extended these results to directed graphs. \color{black}{In  related works, the notion of  zero forcing number of a graph, that is, the minimum cardinality of its ZFS, has been investigated. In particular, it has been shown that by deleting exactly one edge from a graph, 
 the difference between the zero forcing number of the new graph and the old one is between -1 and 1 \cite{edholm2012vertex,berliner2013minimum}.\color{black}  

In this paper, by considering a fixed ZFS as a set of control nodes of a directed network, its ss-controllability is studied. In this regard, we first focus on LTI networks, and  in \S III, we introduce a \emph{$(\mathcal{C},T)$-constructed graph}  and show that there is a one-to-one correspondence between  \color{black} SSC \color{black} networks and $(\mathcal{C},{T})$-constructed graphs. With this in mind, we propose a procedure for synthesizing  networks that are SSC, which is an extension of the results of \cite{mousavi2018structural} to directed networks. 
\color{black} As an extension of the results of \cite{mousavi2017robust}, we also establish a one-to-one correspondence between the set of nodes rendering a family of $(\mathcal{C},T)$-constructed graphs controllable and sources of chains in $\mathcal{C}$.  \color{black}
In the meantime, by considering the time or  iteration by which any color-change force is performed, a method built upon attributing successive integers to the nodes
of \color{black}a SSC \color{black} network is presented that results in a framework for analyzing ss-controllability. 
In particular, our work provides a machinery for studying the robustness of ss-controllability with respect to structural perturbations; see \S IV. 

\color{black}{We note that our approach departs from \cite{edholm2012vertex} and \cite{berliner2013minimum}, since in these works, the effect of removing only one edge of a graph on the zero forcing number has been explored; while,  we consider a fixed  ZFS as the set of control nodes and characterize the maximal set of edges, the addition of 
any subset of which, preserves ss-controllability. }\color{black} 
 In this direction, we introduce the notion of    \emph{perfect graphs}
corresponding to a network with a given set of control nodes. Furthermore, we show that if the number of edges added to a network is greater than a tight bound, ss-controllability of the  network will be compromised. This bound, that is surprisingly independent from the topology of the network, depends only on the size of the network and the number of its control nodes, and it can  increase by enhancing the number of control nodes. Similar results are provided for the removal of edges from a network while ensuring
ss-controllability.

In \S V, we provide an algorithm for combining SSC networks such that the resulting \emph{network-of-networks} is SSC. In~\cite{chapman2014controllability}, controllability of the Cartesian product of networks (with system matrices 
restricted to have a symmetry-preserving property) has been studied. Similarly, in~\cite{yazicioglu2013leader} and \cite{ji2016design}, methods for combining diffusively coupled networks for building a larger controllable network are presented. In this paper, we present a method for combining networks with system matrices of
the same zero/nonzero patterns. Implicit in the proposed method is the existence of structural uncertainties in the 
constituent graphs. In this direction,
we determine the maximum number of edges that can be added between networks while preserving the ss-controllability of the resulting network-of-networks.
 Moreover, as an application of the developed theory, networks with structures described by \emph{directed acyclic graphs} are considered, that appear in hierarchical social networks \cite{lu2015exploring}, and networks with causal inferences \cite{elwert2013graphical}. Based on the properties of directed acyclic graphs, a procedure for combing networks with acyclic graphs is proposed, ensuring the controllability of the overall network with a single control node. 
\color{black}  \S VI is dedicated  to controllability analysis of a family of LTV networks. The work \cite{reissig2014strong} has investigated the controllability of a family of LTV networks whose edge weights vary over time, but the network structure remains unchanged. In \S VI, we generalize the results of \cite{reissig2014strong} to LTV networks that in addition to their edge weights,  have a time-varying structure. \S VII concludes the paper.  \color{black}

\section{Notation and Preliminaries}
In this section, we present the notation and relevant background and constructs 
for our subsequent discussion.

\textit{Notation:} We denote by $\mathbb{R}$ and $\mathbb{Z}$, respectively, the set of real numbers and integer numbers. For a matrix $M\in\mathbb{R}^{p\times q}$, $M_{ij}$ is the entry of $M$ on its $i$th row and $j$th column. The $n\times n$ identity matrix is given by $I_n$, and $e_j$ represents its $j$th column. The cardinality of a set $\mathcal{S}$ is given by $|\mathcal{S}|$. \color{black}For $a,b\in \mathbb{Z}$ with \color{black}$b\geq a$, we define $[a,b]\triangleq\{x\in \mathbb{Z}: a\leq x\leq b\}$. Note that $[a,a]=\{a\}$. \color{black}  Also, for $t_1,t_2\in \mathbb{R}$ with $t_2\geq t_1$, we define $[t_1,t_2]_{\mathbb{R}}=\{x\in\mathbb{R}: t_1\leq x\leq t_2\}$.

 A graph is denoted by $G=(V(G),E(G))$, where $V(G)=\{1,\ldots,n\}$ and $E(G)\subseteq V(G)\times V(G)$ are, respectively, the node set and the edge set of $G$. Also, $|V(G)|$ is the size of $G$. 
 If $(i,j)\in E(G)$, the node $j$ (resp., $i$) is an out-neighbor (resp., in-neighbor) of $i$ (resp., $j$). Any edge $(i,i)\in E(G)$,  $1\leq i\leq n$, is  a self-loop on node $i$. 
Given  graphs $G$ and $G'$, we say that 
 $G\subseteq G' $ if $V(G)=V(G')$ and $E(G)\subseteq E(G')$. For a graph $G$ and a set of edges $E'\subset V(G)\times V(G)$, let $G+E'$ (resp., $G-E'$) be a graph with the node set $V(G)$ and the edge set
 $E(G)\cup E'$ (resp., $E(G)\setminus E'$).

 A chain $C$ is a directed path graph, whose start node (resp., end node) with no in-neighbor (resp., out-neighbor) is called the source (resp., sink). For a node $v\in V(C)$ which is not a sink  of $C$, $v+1$  denotes the out-neighbor  of $v$.

A \color{black}qualitative class \color{black} of  $G$ is a set of patterned matrices defined as $\mathcal{Q}(G) =\{A\in \mathbb{R}^{n\times n}: \mathrm{for}\:\: \: i\neq j,A_{ij}\neq 0\Leftrightarrow(j,i)\in E(G)\}$. 
Note that the diagonal entries of $A\in \mathcal{Q}(G)$ can have any real value.

\subsection{Problem setup}

Consider an LTV network with the following dynamics
\begin{align}
\dot{x}(t)=A(t)x(t)+Bu(t),
\label{lv}
\end{align}
where $x(t)=[x_1(t),\ldots,x_n(t) ]^T$ is the vector of node states, and $u(t)=[u_1(t),\ldots,u_m(t) ]^T$ is the vector of input signals. Moreover, $B\in \mathbb{R}^{n\times m}$ is a constant matrix, called the input matrix,  defined as
$B=[e_{j_1},\ldots,e_{j_m} ]$,
where for $i = 1,\ldots,m$, $j_i \in \{1,...,n\}$. For $i=1,\ldots,m$, the node $j_i$ into which the input signal $i$ is directly injected is called a \textit{control node}. The set of control nodes is represented by $V_C=\{j_1,\ldots,j_m\}$.  In addition, the time-varying matrix $A(t)\in\mathbb{R}^{n\times n}$, for every $t\in\mathbb{R}$,  is a piecewise-continuous function of time, and is called the system matrix. For an LTV system, if $u(t)=0$, for all $t\geq t_0$, one can write the unique solution of (\ref{lv}) as $x(t)=\Phi(t,t_0)x(t_0)$, where the $n\times n$ matrix $\Phi(t,t_0)$ is referred to as  the transition matrix at $(t,t_0)$ and can be written via the Peano-Baker series~\cite{hespanha2018linear} as,
\begin{equation}
\begin{aligned}
\Phi(t,t_0)&=I_n+\int_{t_0}^t A(s_1)ds_1+\int_{t_0}^t A(s_1)\int_{t_0}^{s_1} A(s_2)ds_1ds_2\\
&+\int_{t_0}^t A(s_1)\int_{t_0}^{s_1} A(s_2)\int_{t_0}^{s_2} A(s_3)ds_3ds_2ds_1+\ldots
\end{aligned}
\label{tm}
\end{equation}
When $A(t)$ in (\ref{lv}) is constant, we have an LTI system $
%
\dot{x}(t)=Ax(t)+Bu(t).$
This LTI system, or equivalently the pair $(A,B)$,  is  controllable if there is a suitable  input, steering the states of the system  between any initial and final values within a finite time.

\begin{deff}
Given an LTI system and a graph $G$, the corresponding LTI network (on $G$) is called \color{black}{\textit{strongly structurally controllable} (SSC)} \color{black}if for all $A\in \mathcal{Q}(G)$, the pair $(A,B)$ is controllable.\footnote{\color{black} Note that an LTI network with graph $G$ remains SSC when any self-loop is added to or removed from $G$. }
\label{de1}
\end{deff}

\color{black}

 An LTV system  (\ref{lv}) is said to be controllable on an interval $[t_0,t_1]$ if  there is a suitable input that can derive the system from any initial state at time $t_0$ to any final state at time $t_1$ \cite{sontag2013mathematical,reissig2014strong}. 
%
 The next result presents a controllability condition for an LTV system.
 \begin{pro}[\cite{reissig2014strong}]
Let $t_0,t_1\in \mathbb{R}$, $t_0<t_1$, and $\nu\in\mathbb{R}^n$. The system (\ref{lv}) is controllable on $[t_0,t_1]$ if and only if  $\nu^T \Phi(t_1,\tau)B=0 $, for almost every $\tau\in[t_0,t_1]$, implies that $\nu=0$. 
\label{tm2}
\end{pro}
  
  In this paper, first we consider an LTI network that is  SSC, \color{black} and investigate the preservation of ss-controllability under structural perturbations and combination of networks. \color{black}More precisely, we aim to find sets of edges with the maximum size, adding/deleting any subset of  which to/from a network, does not disturb its ss-controllability property. 
   Moreover, we propose methods for combining SSC networks so that the resulting network is SSC. \color{black} Finally, we propose controllability conditions for a family of LTV networks whose zero-nonzero pattern can also vary over time.
\section{Zero forcing sets}

In this section, zero forcing sets that characterize a set of control nodes, rendering a network SSC, are introduced.
We first review a useful coloring process on the nodes of a graph~\cite{work2008zero}. 

Consider a graph $G$ whose nodes are colored either black or white. \color{black}The color of the nodes can be changed according to a color change-rule:
If a black node $v\in V(G)$ has only one white out-neighbor $u\in V(G)$, it  forces this node to become black; we designate this by $v \rightarrow u$.  The process of repeatedly applying the color-change rule until no more changes are possible is called a forcing process.\color{black}

Now, let $S\subset V(G)$ be the set of initially black nodes of $G$. The set of final black nodes obtained by performing the forcing process is called the derived set of $S$ and is denoted by $\mathcal{D}(S)$.  Given a set of initially black nodes $Z\subset V(G)$, if $\mathcal{D}(Z)=V(G)$, $Z$ is referred to as a zero forcing set (ZFS). 
Considering a forcing process,  a \textit{chronological list of forces}, or simply a \emph{list of forces} $\mathcal{F}$, is defined as a record of the forces in the order in which they are performed. 
 Finally, given a list of forces $\mathcal{F}$, a sequence of nodes $(v_1,\ldots,v_k)$ is a \textit{forcing chain} if for $i=1,\ldots,k-1$,  $v_i\rightarrow v_{i+1}$. 
This forcing chain is \textit{maximal} if $v_1\in Z$, and $v_k$ does not force any node of the graph during the forcing process. \color{black}Note that the maximal forcing chains are node-disjoint in the sense that they do not have any common node. \color{black} In fact, every node of a graph can force at most one other node and can be forced by at most one other node. With this in mind, there are $|Z|$ maximal forcing chains in a graph $G$ (covering all of its nodes)\color{black}\cite{berliner2013minimum}\color{black}. There is a one-to-one correspondence between control nodes rendering a network SSC and the ZFS's.

\begin{theo}[\cite{monshizadeh2014zero}]
An LTI network with the graph $G$ is SSC if and only if $V_C$ is a ZFS of $G$. 
\label{th1}
\end{theo}

\subsection{SSC networks and $(\mathcal{C},T)$-constructed graphs}

Given a set of control nodes, our first contribution is a method for synthesizing directed networks that are SSC. We then present an algorithm for robustness analysis of ss-controllability under edge additions and deletions. First, we review some relevant definitions.

\color{black}Consider a set of \color{black}node-disjoint chains, \color{black} denoted as $\mathcal{C}=\{C_1, \ldots,C_m\}$, where $C_i=(V(C_i), E(C_i))$ is a chain of size $n_i$, $i=1,\ldots,m$. Let $V=\bigcup_{i=1}^m V(C_i)$, and $n=\Sigma_{i=1}^m n_i$. The sources of $\mathcal{C}$ are  the union of sources of $C_i$, $i=1,\ldots,m$. Now, for $\gamma=n-m+1$, let us define the \emph{time function} $T:V\rightarrow [1,\gamma]$ as a function that assigns any node $v\in V$ an integer number $T(v)$, satisfying the following conditions: 1) If for some $C_i$, $i=1,\ldots,m$, $v$ is a source, $T(v)=1$, 2) For any two nodes $u,v\in V$, none of which is a source, one has $T(u)\neq T(v)$, 3) For every $v\in V(C_i)$ that is not a sink ($1\leq i\leq m$), we have $T(v)<T(v+1)$. 

Now, consider a node $v\in V$. Let us define another function $T_{\mathrm{max}}:V\rightarrow [1,\gamma]$ in the following way: If $v$ is a sink of some $C_i$, $i=1,\ldots,m$, define $T_{\mathrm{max}}(v)=\gamma$; otherwise, $T_{\mathrm{max}}(v)=T(v+1)-1$.   

 \color{black}
 For example, the set of chains $\mathcal{C}=\{C_1,C_2\}$ and the time interval $[T(v),T_{\max}(v)]$ for every node $v$ are illustrated in Fig. \ref{design} (a). Then, one can see that $T(v_1)=T(u_1)=1$, $T(v_2)=2$, $T(u_2)=3$, and $T(v_3)=4$. Moreover, $T_{\mathrm{max}}(v_1)=1$, $T_{\mathrm{max}}(u_1)=2$, $T_{\mathrm{max}}(v_2)=3$, $T_{\mathrm{max}}(u_2)=4$, and $T_{\mathrm{max}}(v_3)=4$.
  
 \color{black}
\begin{deff}
For a given set of node-disjoint chains $\mathcal{C}$ and a time function $T$, a \emph{class of $(\mathcal{C}, T)$-constructed graphs}, denoted by $\mathcal{G}^{\mathcal{C},T}$, includes any graph $G$, satisfying the following properties: 1) $V(G)=\bigcup_{i=1}^m V(C_i)$, 2) $\bigcup_{i=1}^m E(C_i)\subseteq E(G)$, 3) For any $u,v\in V(G)$ that $(u,v)\notin \bigcup_{i=1}^m E(C_i)$, we have $(u,v)\notin E(G)$ if $T_{\mathrm{max}}(u)<T(v)$. The sources of a $(\mathcal{C},T)$-constructed graph are the sources of $\mathcal{C}$.
\label{de3}
\end{deff}\color{black}


\begin{figure}[!h]
\centering
\begin{tikzpicture}[scale=.3]
\tikzset{mycirc/.style={draw,shape=circle,style=thick,inner sep=2pt,scale=.7}}
\tikzset{mytext/.style={shape=circle,style=thick}}
    ]
    \node[mycirc,label=above:{\small $[1,1]$},fill=black,text=white] (v1) at (0,0) {\Large $v_1$};
    \node[mycirc,label=above:{\small $[2,3]$}] (v2) at (4,0) {\Large $v_2$};
    \node[mycirc,label=above:{\small $[4,4]$}] (v3) at (8,0) {\Large $v_3$};
\node[mytext] (a) at (4,-8) {(a)};

    \draw[->, line width=2pt] (v1) -- (v2);
    \draw[->, line width=2pt] (v2) -- (v3);

    \node[mycirc,label=below:{\small $[1,2]$},fill=black,text=white] (u1) at (1.5,-4) {\Large $u_1$};
    
    \node[mycirc,label=below:{\small $[3,4]$}] (u2) at (6.5,-4) {\Large $u_2$};

    \draw[->, line width=2pt] (u1) -- (u2);



\begin{scope}[shift={(14.5,0)}]

\tikzset{mycirc/.style={draw,shape=circle,style=thick,inner sep=2pt,scale=.7}}
\tikzset{mytext/.style={shape=circle,style=thick}}
    ]
    \node[mycirc,label=above:{\small $[1,1]$},fill=black,text=white] (v1) at (0,0) {\Large $v_1$};
   \node[mycirc] (v2) at (4,0) {\Large $v_2$};
    \node[mycirc,label=above:{\small $[4,4]$}] (v3) at (8,0) {\Large $v_3$};
\node[mytext] (b) at (4,-8) {(b)}; 

\node[mytext] (c) at (4,2.5) {$[2,3]$};

    \draw[->, line width=2pt] (v1) -- (v2);
    \draw[->, line width=2pt] (v2) -- (v3);

    \node[mycirc,label=below:{\small $[1,2]$},fill=black,text=white] (u1) at (1.5,-4) {\Large $u_1$};
    
    \node[mycirc,label=below:{\small $[3,4]$}] (u2) at (6.5,-4) {\Large $u_2$};

    \draw[->, line width=2pt] (u1) -- (u2);
    
     \draw[<->, dashed, style=thick] (u1) -- (v1);
      \draw[<->, dashed, style=thick] (u1) -- (v2);
       \draw[->, dashed, style=thick] (v3) -- (u1);
    \draw[->, dashed, style=thick] (u2) -- (v1);
     \draw[<->, dashed, style=thick] (u2) -- (v2);
      \draw[<->, dashed, style=thick] (u2) -- (v3);

      \draw[->,dashed,line width=1pt] (u1) to [out=235,in=185,looseness=8] (u1);
      \draw[->,dashed,line width=1pt] (v1) to [out=125,in=175,looseness=8] (v1);
      \draw[->,dashed,line width=1pt] (v3) to [out=55,in=5,looseness=8] (v3);
      \draw[->,dashed,line width=1pt] (v2) to [out=65,in=115,looseness=7] (v2);
      \draw[->,dashed,line width=1pt] (u2) to [out=305,in=355,looseness=8] (u2);

      \draw[dashed,line width=1pt,->,shorten >=1pt] (v2) to [out=150,in=30] (v1);
      \draw[dashed,line width=1pt,->,shorten >=1pt] (v3) to [out=150,in=30] (v2); 
      \draw[dashed,line width=1pt,->,shorten >=1pt] (u2) to [out=210,in=330] (u1); 
       \draw[dashed,line width=1pt,->,shorten >=1pt] (v3) to [out=135,in=45] (v1);

\end{scope}
\end{tikzpicture}
\caption{a) Set of node-disjoint chains $\mathcal{C}$ and time interval $[T(v),T_{\max}(v)]$, b) $(\mathcal{C},T)$-constructed graph.}
\label{design}
\end{figure}
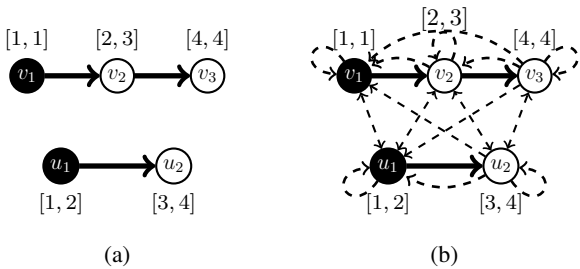


In Fig.~\ref{design} (b),  the solid directed lines, that is, the union of edges of $C_1$ and $C_2$, denote edges which \textit{should} exist in a $(\mathcal{C},T)$-constructed graph. Moreover, dotted lines show edges that \textit{can} exist in this graph.

 The next result demonstrates that by synthesizing a $(\mathcal{C},T)$-constructed graph and choosing its sources as control nodes, \color{black}a SSC \color{black} network is provided.
 
 \begin{theo}
The set of sources of a $(\mathcal{C},T)$-constructed graph is \color{black}a ZFS of this graph. \color{black}
\label{th2}
 \end{theo}

\textit{Proof:} Let $G\in\mathcal{G}^{\mathcal{C},{T}}$, and $Z$ be the set of its sources, but suppose $Z$ is not a ZFS, i.e., $\mathcal{D}(Z)\neq V(G)$. 
 For $i=1,\ldots,m$, let $u^i_L$ be a node in $C_i\cap \mathcal{D}(Z) $ such that $u^i_L+1$ is white. In other words, $u^i_L$ is the terminating node in $C_i$ that is colored black. Let $\theta^i_L=T_{\mathrm{max}}(u^i_L)$. Without loss of generality, assume that $\theta^1_L=\mathrm{min} _{1\leq i\leq m} \theta^i_L $. Since $u^1_L$ does not force any node to become black, it has at least two white out-neighbors, one of which is in $C_1$. Let $s\in C_j$,  for some $1<j\leq m$, be its other white out-neighbor. Then, $T(s)>\theta^j_L\geq \theta^1_L$. Thus, $T_{\mathrm{max}}(u^1_L)< T(s)$, and based on Definition \ref{de3}, $(u^1_L,s)\notin E(G)$, which is a contradiction. \carre

Now, in a reverse direction, we show that  every  SSC network corresponds to a $(\mathcal{C}, T)$-constructed graph. Consider  a chronological list of forces $\mathcal{F}$ in \color{black}an LTI SSC \color{black} network. \color{black}We will show that the robustness analysis of ss-controllability can be facilitated by considering the time or the iteration index in which any node is forced. Note that in every iteration of the forcing process, there may be more than one potential force that can be independently performed at the same time. For example, in \cite{hogben2012propagation,kenter2018error}, the propagation time of a ZFS, that is, the minimum number of iterations for simultaneous application of the color-change rule until the termination of the forcing process has been studied. \color{black}However, to the purpose of this paper, we allow only one node to be colored black at each step of the forcing process. \color{black}In this direction, Algorithm 1 shows how every node $v$ of \color{black}a SSC \color{black} network can be assigned a \emph{forcing time} $T(v)$. Indeed, $T(v)$ is the iteration index  in which $v$ becomes black. Note that the computational complexity of Algorithm 1 is $\mathcal{O}(n^2)$, since it is the same as the  computational complexity of coloring a graph through a forcing process \cite{brimkov2019computational}. \color{black} 

 \begin{figure}[htb]
\mbox{}\hrulefill
\vspace{-.4em}
\\
\textbf{Algorithm 1:}

\vspace{-.7em}
\mbox{}\hrulefill
\vspace{-.01em}
\\
\small Given $Z$ \\
\begin{tabular}{l}
For every $v\in Z$\\
$T(v)=1$;\\
$\gamma=1$;\\
\end{tabular}
\\
while the color change rule is possible, do\\
\begin{tabular}{|l}
Among nodes that can force their out-neighbors,\\ choose exactly one node $w$, and let $w\rightarrow u$.\\ 
$T(u)=\gamma+1$.\\
%
$\gamma=\gamma+1$;\\
\end{tabular}\\
end while\\
\vspace{-.5em}
\mbox{}\hrulefill
\caption{An algorithm that endows nodes of a graph with  successive integers with respect to a given ZFS.}
\vspace{-0.08in}
\label{A1}
\end{figure}

\begin{theo}
Consider an LTI SSC \color{black} network with the graph $G$ and the set of control nodes $Z$. Then, there is  a time function $T:V(G)\rightarrow [1,\gamma]$  and  a set of node-disjoint chains $\mathcal{C}$ with the sources $Z$ such that $G\in \mathcal{G}^{\mathcal{C},T}$.
\label{theorem_3}  
\end{theo}

\emph{Proof:} \color{black}Consider a list of forces $\mathcal{F}$ associated with $Z$, and let $|Z|=m$. Define $\mathcal{C}=\{C_1,\ldots,C_m\}$ as the set of maximal forcing chains associated with $\mathcal{F}$. \color{black} Also, let $T(.)$ be a function that assigns every node $v$ a forcing time $T(v)$, as provided by Algorithm 1.  Now, suppose that $G\notin \mathcal{G}^{\mathcal{C},T}$. Then, either (1) for some $i\neq j$, there are some $u\in V(C_i)$ and $v\in V(C_j)$ such that $T_{\mathrm{max}}(u)<T (v)$, and $(u,v)\in E(G)$, or (2) there are some $u, v\in V(C_k)$ such that $T_{\mathrm{max}}(u)<T (v)$ and $(u,v)\notin E(C_k)$, but $(u,v)\in E(G)$. In both cases, note that $u$ is the last black node of $C_i$ in time $T_{\mathrm{max}}(u)$, since the node $u+1$ is forced in the time $T_{\mathrm{max}}(u)+1$. 
 However, $u$ has another white out-neighbor $v$ which will become black in time $T(v)$, and note that $T_{\mathrm{max}}(u)<T(v)$. Hence, $u$ has two white out-neighbors in time $T_{\mathrm{max}}(u)$, and the color change rule cannot be performed in this step, establishing a contradiction. \color{black} \carre

\color{black}For an LTI SSC \color{black} network with the graph $G$, consider a list of forces $\mathcal{F}$, and  let $\mathcal{C}$ and $T$ be, respectively, a set of node-disjoint chains and a time function associated with $\mathcal{F}$ such that $G\in\mathcal{G}^{\mathcal{C},T}$. Then, for every $v\in V(G)$, $T(v)$ is the time step in which the node $v$ is colored black. Let for some $1\leq i\leq m$, $v\in C_i$. Then, $T_{\mathrm{max}}(v)+1$ is the time step when $v $ forces $v+1\in C_i$ to be black. One can see that  node $v$ is the last black node of the chain $C_i$ during the time interval $\mathbb{T}(v)=[T(v),T_{\mathrm{max}}(v)]$. Accordingly, we can have $(u,v)$ as an edge of the network without disturbing the ss-controllability if $v$ is forced before $u$ forces, or equivalently, if $T_{\mathrm{max}}(u)\geq T(v)$.
For any $u,v\in V(G)$ for which $\mathbb{T}(u)\cap \mathbb{T}(v)\neq \emptyset$, $(u,v)$ and $(v,u)$ can simultaneously exist in \color{black}a SSC \color{black} network. Otherwise, either $T_{\mathrm{max}}(u)> T(v)$ or $T_{\mathrm{max}}(v)> T(u)$. Thereby, For any $(u,v)\in V(G)\times V(G)$, if the network with structure $G+\{(u,v)\}$ is no longer SSC, then a network with the graph $G+\{(v,u)\}$ would be SSC. 

\color{black} Next, we define a new family of matrices, associated with a class of $(\mathcal{C},T)$-constructed graphs, whose zero-nonzero pattern is not necessarily the same.

\begin{deff}
A qualitative class corresponding to $\mathcal{G}^{\mathcal{C},T}$ is defined  as $\mathcal{P}(\mathcal{G}^{\mathcal{C},T})=\{A\in \mathbb{R}^{n\times n}: A\in \mathcal{Q}(G), \quad \mathrm{for}\quad\mathrm{some}\quad G\in\mathcal{G}^{\mathcal{C},T}\}$.
\end{deff}

The next theorem provides a necessary and sufficient controllability condition for every LTI network whose system matrix is in $\mathcal{P}(\mathcal{G}^{\mathcal{C},T})$.
\begin{theo}
Any LTI system  with $A\in\mathcal{P}(\mathcal{G}^{\mathcal{C},T})$ is controllable if and only if $V_C$ is any set that includes the set of sources of $\mathcal{C}$.
\label{p}
\end{theo}

\emph{Proof:} Let $\mathcal{S}$ be the set of sources of $\mathcal{C}$. Consider a graph $G$, where $E(G)=\bigcup_{i=1}^m E(C_i)$, and note that $G\in \mathcal{G}^{\mathcal{C},T}$. From Theorem \ref{th1}, if an LTI system  with  $A\in\mathcal{Q}(G)$  is SSC,  $V_C$ is a ZFS of $G$. Then, since  any  ZFS of $G$ should include $\mathcal{S} $, the necessary condition is proved. Moreover, since $\mathcal{S}$ is a ZFS of any $G\in\mathcal{G}^{\mathcal{C},T}$, Theorem \ref{th1} implies that for $V_C=\mathcal{S}$, any LTI system  with $A\in \mathcal{Q}(G)$ is controllable. 
\carre

\color{black} In the next section, we study the maximum number of edges that can be added to a network while the ss-controllability is preserved.\color{black} 

\section{Robustness of SS-Controllability }
Strong structural controllability captures network controllability with implicit robustness guarantees  against variations in the edge weights. In this context, however, no edges are allowed to be added to, or removed from the network. We address this shortcoming next.
In this regard, we introduce the notion of a critical additive (resp., subtractive) edge-set,  \color{black} a  set of edges of the maximum cardinality \color{black} whose any subset can be added to (resp., removed from) a network with a given set of control nodes, while the ss-controllability is  preserved. Some intermediate notions are first introduced.

\color{black}For a graph $G$ with a given ZFS, consider a list of forces $\mathcal{F}$ and the associated set  of node-disjoint chains $\mathcal{C}$. Then, provide a time function $T:V\rightarrow [1,\gamma]$ according to Algorithm 1 in Fig \ref{A1}.
\color{black}


\begin{deff}
\color{black}A $(\mathcal{C},T)$-constructed graph $G$  
is called perfect if  $T_{\mathrm{max}}(u)\geq T(v)$ implies that $(u,v)\in E(G)$, for any $u,v\in V(G)$. 
We denote a perfect $(\mathcal{C},T)$-constructed graph and its edge set respectively  by $\mathcal{G}^{\mathcal{C},T}_{\mathrm{perf}}$ and  $E_{\mathrm{perf}}$.
\label{de4}
\end{deff}
\color{black}

For example, the graph in Fig.~\ref{design} (b), with all dotted and solid lines as its edges, represents a perfect $(\mathcal{C},T)$-constructed graph.

\color{black}\begin{rem}
Note that with  a set of node-disjoint chains $\mathcal{C}$ and a time function $T$, only a unique $\mathcal{G}^{\mathcal{C},T}_{\mathrm{perf}}$ can be defined.  \color{black} However,  
  if $|\cup_{i=1}^m E(C_i)|=q$, and $|E_{\mathrm{perf}}|=l$, then there are $2^{l-q}$  graphs, all of which are $(\mathcal{C},T)$-constructed.\color{black}
\end{rem}  \color{black}

\begin{pro}
If $G=\mathcal{G}^{\mathcal{C},T}_{\mathrm{perf}}$, \color{black}then for every set of node-disjoint chains $\mathcal{C}'$ and any time function $T'$ that $G\in \mathcal{G}^{\mathcal{C}', T'}$\color{black}, we have $T'(v)=T(v)$, for all $v\in V(G)$.
\label{pr1}
\end{pro}

\textit{Proof:}  Since $G=\mathcal{G}^{\mathcal{C},T}_{\mathrm{perf}}$, one of its associated time function is $T$.   Now, we claim that in every iteration $k$, $1< k\leq \gamma$, there is exactly one white node that can be forced to be black. The proof follows by contradiction. Assume that the iteration $i$, $1<i\leq \gamma$, is the first iteration in which there
are at least two different white nodes $v_1$ and $v_2$ that can be forced by two black nodes $u_1$ and $u_2$. Note that we should have $u_1\neq u_2$; since any node can force at most one white node.   Now, let $v_1\in V(G)$ be a node \color{black}such \color{black} that $T(v_1)=i$.   Since $T(v_1)=i$, $u_2$  remains as the last black node of its chain in the iteration $i$, and then $T_{\max}(u_2)\geq i$. Thus, $(u_2,v_1)\in E(G)$, implying that $u_2$ has two white out-neighbors $v_1$ and $v_2$ in the iteration $i$; hence, it cannot force $v_2$ in this iteration, which contradicts the assumption.  \carre 

A graph $G$ with a given ZFS is said to have a perfect graph $G'$, if for a set of node-disjoint chains $\mathcal{C}$ with the set of sources ZFS and a time function $T$, $G'=\mathcal{G}^{\mathcal{C},T}_{\mathrm{perf}}$, and $G$ is a $(\mathcal{C},T)$-constructed graph. Note that since in every step of Algorithm 1, exactly one node is \emph{chosen} from among the nodes having one white out-neighbor, the resulting set of chains $\mathcal{C}$ and the time function $T$ are not unique in general. Then, a graph may have different perfect graphs. However, with the aid of Proposition \ref{pr1}, one can see that if for some $\mathcal{C}$ and $T$, $G=\mathcal{G}^{\mathcal{C},T}_{\mathrm{perf}}$, then it has a unique perfect graph.


 \color{black} Assume that there is an edge, which when added to a SSC network, the ss-controllability is preserved. Next, we show that this edge  belongs to the edge set of one of the associated perfect graphs.
 \color{black}  


\begin{lem}
Consider an LTI SSC \color{black} network with the graph $G$. Let  $(u,v)\in V(G)\times V(G)$, and $(u,v)\notin E(G)$. If the network with the graph $G'=G+\{(u,v)\}$ and control nodes $V_C$ is SSC, \color{black}then one can find a time function $T$ and a set of node-disjoint chains $\mathcal{C}$ with the set of sources $V_C$, \color{black} such that $G$ and $G'$ are both $(\mathcal{C},T)$-constructed graphs.
\label{pr2}
\end{lem}

\textit{Proof:} \color{black} It suffices to show that there is a list of forces $\mathcal{F}$ for both $G$ and $G'$, in which $u\nrightarrow v$. In other words, we should prove that $(u,v)\notin  E(\mathcal{C})$, where $\mathcal{C}$ is the set of maximal forcing chains associated with $\mathcal{F}$. Consider a list of forces $\mathcal{F}^*$ with the set of node-disjoint chains $\mathcal{C}^*=\{C^*_1,\ldots,C_m^*\}$ and the time function $T^*$ such that $G\in \mathcal{G}^{\mathcal{C}^*,T^*}$.  Apply the forces of $\mathcal{F}^*$ in $G'$ until the iteration $T_{\mathrm{max}}^*(u)$, and let $\mathcal{B}_1$ be the set of black nodes until this iteration.  If $v$ is forced before $u$ forces, then $u\rightarrow v$ is not in the list of forces $\mathcal{F}^*$, and thus, one can define $\mathcal{F}=\mathcal{F}^*$. Now, assume that in the iteration $T^*_{\mathrm{max}}(u)$, $v$ is white, that is, $T^*_{\mathrm{max}}(u)< T^*(v)$.   Let $u\in C^*_k$, for some $1\leq k\leq m$. \color{black}Note that $u$ is not the sink of $C^*_k$; otherwise, $T^*_{\mathrm{max}}(u)= \gamma$ and  $T^*_{\mathrm{max}}(u)\geq T^*(z)$, for every $z\in V(G)$. 
 \color{black} Now, suppose that $u$ has another white out-neighbor in $C_k^*$, say $w$. Then, in the iteration $T^*_{\mathrm{max}}(u)$, $u$ has two white out-neighbors $v$ and $w$ in $G'$. Thus, at least one of these neighbors should be forced by another node of the graph. If $v$ is forced by another node except $u$, then a list of forces can be found in which the force $u\rightarrow v$ does not appear, and the proof is complete. Now, suppose that $w$ is forced by another node of $G'$ other than $u$. 
\color{black}
Then, after the iteration $T_{\max}(u)$,  a subsequence of forces as $\mathcal{F}_1=(z_1\rightarrow z_2, \ldots, z_q\rightarrow w)$ can be found in $G'$, which causes $w$ to become black. Note that none of $z_i$'s, $1\leq i \leq q$, equals to $u$ or $v$.
\color{black}
Now, consider a list of forces in $G$ after the iteration in which $u$ forces $w$ until the iteration in which $v$ becomes black. In this direction, for $y_i\in V(G)$, $1\leq i\leq p$, let $ \mathcal{F}^*_1= (y_1\rightarrow y_2,y_2\rightarrow y_3,\ldots, y_{p}\rightarrow v )$, which includes the forces  from the iteration $T^*(w)+1$ until the iteration $T^*(v)$. Then, none of $y_i$'s is the same as the node $u$.  Now, we claim that after $w$ becomes black in $G'$, the sequence of forces $ \mathcal{F}^*_1$ can be performed in $G'$ as well. Assume that the claim is not true, and let $k$, $1\leq k\leq p$, be the smallest index, where the force $y_k\rightarrow y_{k+1}$ cannot be performed in $G'$.  It   implies  that $y_k$ has at least two white out-neighbors in $G'$. Note that the only node in $G'$ that has a new out-neighbor is $u$, and  $y_k\neq u$. Let $\mathcal{B}_2$ and $\mathcal{B}'_2$ be the set of black nodes respectively in $G$ and $G'$ before applying the force $y_k\rightarrow y_{k+1}$. Then, based on the previous discussion,  one can see that  $\mathcal{B}_2=\mathcal{B}_1\cup\{w\}\cup \{y_1,\ldots,y_k\}$. In other words, $\mathcal{B}_2$ includes the nodes in $G$ that become black before $u$ forces $w$, together with  node $w$, and the set of nodes $y_i$'s, $1\leq i\leq k$, that become black  before $y_k$ forces $y_{k+1}$. Moreover, we have $\mathcal{B}'_2=\mathcal{B}_2\cup\{z_1,\ldots,z_q\}$. In fact, one can see that $\mathcal{B}'_2$ includes all nodes of $\mathcal{B}_1$,  all the nodes $z_i$'s, $1\leq i\leq q$, that become black before $w$ is forced, and  node $w$. Also, since based on the assumption, for all $1\leq i< k$, the force $y_i\rightarrow y_{i+1}$ can be performed in $G'$, $\mathcal{B}'_2$ includes the set of nodes $\{y_1, \ldots, y_k\}$ as well. Thus, $\mathcal{B}_2\subseteq \mathcal{B}'_2$, which implies that
 all  nodes being  black in $G$ before applying the force $y_k\rightarrow y_{k+1}$  are black in $G'$ as well.
 Therefore, the force $y_k\rightarrow y_{k+1}$ can be performed in $G'$ as it can be performed in $G$,  contradicting the assumption. Hence, when $w$ is black, $v $ can be forced by a node other than $u$, completing the proof. \color{black} \carre
%
%

The next theorem is one of our main results and describes networks that do not remain SSC under any new edge addition.
\begin{theo}
Consider an LTI network with the graph $G$. Let $V_C$ be a ZFS of $G$.  By adding any single edge to $G$, the new network is no longer SSC if and only if 
for a time function $T$ and a set of node-disjoint chains $\mathcal{C}$ with sources $V_C$,
$G=\mathcal{G}^{\mathcal{C},T}_{\mathrm{perf}}$.
\label{th3} 
\end{theo}

\textit{Proof:} \color{black}Suppose that 
$G$ is a perfect graph, but for some $(u,v)\in V(G)\times V(G)$ for which $(u,v)\notin E(G)$, a network with the same set of control nodes and the structure $G'=G+\{(u,v)\}$ is SSC. Then from Lemma \ref{pr2}, one can find a time function $T$ and a set of node-disjoint chains $\mathcal{C}$ such that both $G$ and $G'$ are $(\mathcal{C},T)$-constructed. Hence, we have $T_{\mathrm{max}}(u)\geq T(v)$.  Moreover, from Proposition \ref{pr1}, one can see that since $G$ is a perfect graph, for any time function $T'$ and a set of chains $\mathcal{C}'$ that $G\in \mathcal{G}^{\mathcal{C}', T'}$, one has $T'(w)=T(w)$, 
$\forall w\in V(G)$. Thus, since $T_{\mathrm{max}}(u)\geq T(v)$, we have $(u,v)\in E(G)$,  contradicting the assumption.  \color{black}
Now, assume by adding any single edge to $G$, the new network is no longer SSC, but there is not any $\mathcal{C}$ and $T$ for which $G$ is a perfect $(\mathcal{C},T)$-constructed graph. From Theorem \ref{theorem_3}, there is a set of chains $\mathcal{C}'$ and a time function  $T'$ such that $G\in\mathcal{G}^{\mathcal{C}',T'}$. Let $G^*=\mathcal{G}^{\mathcal{C}',T'}_{\mathrm{perf}}$. Since $G\neq G^*$, $G\subset G^*$. Then, $E_{\mathrm{dif}}=E(G^*)\setminus E(G)\neq \emptyset$. Note that a network with the graph $G^*$ and the control nodes $V_C$ is SSC. Hence, for any $e\in E_{\mathrm{dif}}$, the network with the graph $G+\{e\}$ would be SSC, establishing a contradiction.      
\carre

Now, consider \color{black}a SSC \color{black} network with the graph $G$ of size $n$ and a set of control nodes of size $m$. Next, we show that although $G$ may not have a unique perfect graph, the cardinality of the edge set of all of its perfect graphs is the same.   

\color{black}\begin{lem}
\color{black} Consider a graph $G$ of size $n$ with a ZFS of size $m$. Any perfect $(\mathcal{C},T)$-constructed graph of $G$ has an edge set of size \color{black}$|E_{\mathrm{perf}}|=\frac{1}{2} n(n+1)+\frac{1}{2} m(2n-m-1)$\color{black}. 
\label{th4}
\end{lem}

\color{black}\textit{Proof:}  One can partition  $E_{\mathrm{perf}}$ into two sets of edges $E_1$, $E_2$, where $E_1=\bigcup_{i=1}^m E(C_i)$, and $E_2=E_{\mathrm{perf}}\setminus E_1$. We have $|E_1|=n-m$. In addition,  $(u,v)\in E_2$ if and only if $T_{\max}(u)\geq T(v)$. Note that if $u$ is a sink of some $C_i\in\mathcal{C}$, then $T_{\max}(u)=n-m+1$, and there are $m$ chains with $m$ sinks. Then, for any $u$ that is a sink, we have $(u,v)\in E_{\mathrm{perf}}$, for all $v\in  V(G)$. Now, consider a node $u$ that is not a sink. Then,  $T_{\max}(u)=T(u+1)-1$, where $u+1$ is the out-neighbor of $u$ in the same chain. Thus, for any $1\leq k\leq n-m$, there is only a single node $u$ that $T_{\max}(u)=k$.  Let $V^*_u=\{v\in V(G): (u,v)\in E_{\mathrm{perf}}\} $. Then, $v\in V^*_u$ if and only if $T(v)\leq k $. For any node $v$ in the set of sources of $\mathcal{C}$, we have $T(v)=1$. Moreover, for any $1<j\leq k$, there is only one node $v'$ that $T(v')=j$. Thus, $|V^*_u|=m+k-1$. Accordingly, one can write $|E_2|=mn+\sum_{k=1}^{n-m}m+k-1$. Hence, we have $|E_{\mathrm{perf}}|=|E_1|+|E_2|=\frac{1}{2} n(n+1)+\frac{1}{2} m(2n-m-1)$. \carre

It is deduced from Lemma \ref{th4} that for all sets of chains $\mathcal{C}$ and time functions $T$, the number of edges of a perfect $(\mathcal{C},T)$-constructed graph depends only on the size of the graph and the cardinality of its ZFS, independent from the choice of $\mathcal{C}$ and $T$. 
We note that with a single control node, $\frac{|E_{\mathrm{perf}}|}{n^2}$ converges to $0.5$ as $n$ increases. Moreover, by increasing the number of control nodes, 
more edges are available in the corresponding perfect graph.
%
 
\color{black} Lemma  \ref{th4} also provides a method to check that whether a graph of size $n$ and a given ZFS of size $m$ is a perfect $(\mathcal{C},T)$-constructed graph or not. In fact, if $G\in \mathcal{G}^{\mathcal{C},T}$, but $G\neq \mathcal{G}^{\mathcal{C},T}_{\mathrm{perf}}$, then $G\subset\mathcal{G}^{\mathcal{C},T}_{\mathrm{perf}}$, and accordingly $|E(G)|<|E_{\mathrm{perf}}|$. 

\begin{corr}
Consider a graph $G$ of size $n$ and a ZFS of size $m$. For a time function $T$ and a set of node-disjoint chains $\mathcal{C}$, $G=\mathcal{G}^{\mathcal{C},T}_{\mathrm{perf}}$ if and only if $|E(G)|=|E_{\mathrm{perf}}|$. 
\end{corr}
 
 \color{black}

Now, consider the next definition in an LTI network with  
 graph $G$.

\begin{deff} 
Consider an LTI SSC \color{black} network with a set of control nodes $V_C$. A set of edges, adding (resp., removing) any subset of which to (resp., from) the network preserves ss-controllability, is referred to as an \emph{additive edge-set} (resp., \emph{subtractive edge-set}). 
In other words, a set of edges $E^*\subset V(G)\times V(G)$  is an additive (resp., subtractive) edge-set of $G$ if for any $E'\subseteq E^*$, the network with the graph $G+E'$ (resp., $G-E'$) is SSC. The \emph{critical additive number} (resp., \emph{critical subtractive number}) is the maximum of $|E^*|$ over all additive (resp., subtractive) edge-sets $E^*\subset V(G)\times V(G)$ and is denoted by $n^c_{\mathrm{add}}(G)$ (resp., $n^c_{\mathrm{sub}}(G)$). \color{black}A critical additive (resp., subtractive) edge-set is an additive (resp., subtractive) edge-set of the maximum cardinality and is represented by  $E^c_{\mathrm{add}}(G)$ (resp., $E^c_{\mathrm{sub}}(G)$). \color{black} Thus, we have $|E^c_{\mathrm{add}}(G)|=n^c_{\mathrm{add}}(G)$ (resp., $|E^c_{\mathrm{sub}}(G)|=n^c_{\mathrm{sub}}(G)$).  
\label{de5}
\end{deff} 

\begin{theo}
\color{black}For an LTI SSC \color{black} network with the graph $G$, consider a  set of node-disjoint chains $\mathcal{C}$ with the set of sources $V_C$ and a time function $T$ such that $G\in\mathcal{G}^{\mathcal{C},T}$. Then, $E(\mathcal{G}^{\mathcal{C},T}_{\mathrm{perf}})\setminus  E(G)$ is a critical additive edge-set. \color{black} Moreover, \color{black}$n^c_{\mathrm{add}}(G)=\frac{1}{2} n(n+1)+\frac{1}{2} m(2n-m-1) - |E(G)|$\color{black}. 
\label{th5}
\end{theo}

\textit{Proof:} Note that $\bigcup_{i=1}^m E(C_i)\subseteq E(G)\subseteq E_{\mathrm{perf}}$, and let $E^*=E_{\mathrm{perf}}\setminus E(G)$. Then, by Definitions \ref{de3} and \ref{de4}, $E^*$ includes edges which may or may not exist in a $(\mathcal{C},T)$-constructed graph. In other words, by adding any subset of $E^*$ to the edge set of $G$, a new $(\mathcal{C},T)$-constructed graph is obtained which according to Theorems 1 and \ref{th2} is SSC. 
\color{black}Moreover, from Theorem \ref{th2}, by adding any new edge to a network, it would no longer be SSC  if and only if it is a perfect graph. Lemma \ref{th4}, then, implies that any network with a perfect graph of size $n$ and the set of control nodes of size $m$ has an edge set of size $|E_{\mathrm{perf}}|=\frac{1}{2} n(n+1)+\frac{1}{2} m(2n-m-1)$. Thus, the maximum number of edges that can be added to a network without disturbing its ss-controllability is $|E_{\mathrm{perf}}|-|E(G)|$. Now, since $|E^*|=|E_{\mathrm{perf}}|-|E(G)|$, the proof \color{black}is complete\color{black}. \carre   

 Note that by adding any set of edges $E'$ to $G$ with $|E'|>n^c_{\mathrm{add}}(G)$, \color{black} the network would no longer remain SSC. \color{black}

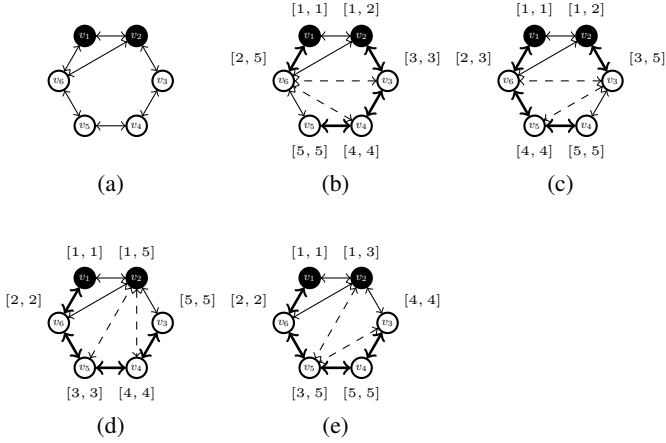
\begin{figure}[!h]
\centering
\begin{tikzpicture}[scale=.23]
\tikzset{mycirc/.style={draw,shape=circle,style=thick,inner sep=2pt,scale=.5}}
\tikzset{mytext/.style={shape=circle,style=thick}}
    ]
    
    \node (v3)[mycirc] at ( 0:3) {$v_3$};
\node[mycirc,fill=black,text=white] (v2) at ( 60:3) {$v_2$};
\node[mycirc,fill=black,text=white] (v1) at (2*60:3) {$v_1$};
\node[mycirc] (v6) at (3*60:3) {$v_6$};
\node[mycirc] (v5) at (4*60:3) {$v_5$};
\node[mycirc] (v4) at (5*60:3) {$v_4$};

\node[mytext] (a) at (270:6) {(a)};

    \draw[<->] (v1) -- (v2);
    \draw[<->] (v2) -- (v3);
    \draw[<->] (v4) -- (v3);
    \draw[<->] (v4) -- (v5);
    \draw[<->] (v5) -- (v6);
    \draw[<->] (v2) -- (v6);
    \draw[<->] (v1) -- (v6);



\begin{scope}[shift={(13,0)}]

\tikzset{mycirc/.style={draw,shape=circle,style=thick,inner sep=2pt,scale=.5}}
\tikzset{mytext/.style={shape=circle,style=thick}}
    ]

    \node (v3)[mycirc,label=above right:{\tiny $[3,3]$}] at ( 0:3) {$v_3$};
\node[mycirc,label=above:{\tiny $[1,2]$},fill=black,text=white] (v2) at ( 60:3) {$v_2$};
\node[mycirc,label=above:{\tiny $[1,1]$},fill=black,text=white] (v1) at (2*60:3) {$v_1$};
\node[mycirc,label=above left:{\tiny $[2,5]$}] (v6) at (3*60:3) {$v_6$};
\node[mycirc,label=below:{\tiny $[5,5]$}] (v5) at (4*60:3) {$v_5$};
\node[mycirc,label=below:{\tiny $[4,4]$}] (v4) at (5*60:3) {$v_4$};

\node[mytext] (b) at (270:6) {(b)};

    \draw[<->] (v1) -- (v2);
    \draw[<->, line width=1pt] (v2) -- (v3);
    \draw[<->, line width=1pt] (v4) -- (v3);
    \draw[<->, line width=1pt] (v4) -- (v5);
    \draw[<->] (v5) -- (v6);
    \draw[<->] (v2) -- (v6);
    \draw[<->, line width=1pt] (v1) -- (v6);
    
     \draw[<->,dashed] (v3) -- (v6);
       \draw[<->,dashed] (v4) -- (v6);
\end{scope}

\begin{scope}[shift={(26,0)}]

\tikzset{mycirc/.style={draw,shape=circle,style=thick,inner sep=2pt,scale=.5}}
\tikzset{mytext/.style={shape=circle,style=thick}}
    ]

    \node (v3)[mycirc,label=above right:{\tiny $[3,5]$}] at ( 0:3) {$v_3$};
\node[mycirc,label=above:{\tiny $[1,2]$},fill=black,text=white] (v2) at ( 60:3) {$v_2$};
\node[mycirc,label=above:{\tiny $[1,1]$},fill=black,text=white] (v1) at (2*60:3) {$v_1$};
\node[mycirc,label=above left:{\tiny $[2,3]$}] (v6) at (3*60:3) {$v_6$};
\node[mycirc,label=below:{\tiny $[4,4]$}] (v5) at (4*60:3) {$v_5$};
\node[mycirc,label=below:{\tiny $[5,5]$}] (v4) at (5*60:3) {$v_4$};

\node[mytext] (c) at (270:6) {(c)};

    \draw[<-> ] (v1) -- (v2);
    \draw[<->, line width=1.pt] (v1) -- (v6);
    \draw[<->, line width=1.pt] (v6) -- (v5);
    \draw[<->, line width=1.pt] (v4) -- (v5);
    \draw[<-> ] (v3) -- (v4);
    \draw[<->] (v2) -- (v6);
    \draw[<->, line width=1.pt] (v2) -- (v3);
    
     \draw[<->,dashed] (v3) -- (v6);
       \draw[<->,dashed] (v3) -- (v5);
\end{scope}

\begin{scope}[shift={(0,-14)}]

\tikzset{mycirc/.style={draw,shape=circle,style=thick,inner sep=2pt,scale=.5}}
\tikzset{mytext/.style={shape=circle,style=thick}}
    ]

    \node (v3)[mycirc,label=above right:{\tiny $[5,5]$}] at ( 0:3) {$v_3$};
\node[mycirc,label=above:{\tiny $[1,5]$},fill=black,text=white] (v2) at ( 60:3) {$v_2$};
\node[mycirc,label=above:{\tiny $[1,1]$},fill=black,text=white] (v1) at (2*60:3) {$v_1$};
\node[mycirc,label=above left:{\tiny $[2,2]$}] (v6) at (3*60:3) {$v_6$};
\node[mycirc,label=below:{\tiny $[3,3]$}] (v5) at (4*60:3) {$v_5$};
\node[mycirc,label=below:{\tiny $[4,4]$}] (v4) at (5*60:3) {$v_4$};

\node[mytext] (d) at (270:6) {(d)};

    \draw[<->] (v1) -- (v2);
    \draw[<->, line width=1pt] (v1) -- (v6);
    \draw[<->, line width=1pt] (v6) -- (v5);
    \draw[<->, line width=1pt] (v4) -- (v5);
    \draw[<->] (v3) -- (v2);
    \draw[<->] (v2) -- (v6);
    \draw[<->, line width=1pt] (v4) -- (v3);
    
     \draw[<->,dashed] (v5) -- (v2);
       \draw[<->,dashed] (v2) -- (v4);
\end{scope}

\begin{scope}[shift={(13,-14)}]

\tikzset{mycirc/.style={draw,shape=circle,style=thick,inner sep=2pt,scale=.5}}
\tikzset{mytext/.style={shape=circle,style=thick}}
    ]

    \node (v3)[mycirc,label=above right:{\tiny $[4,4]$}] at ( 0:3) {$v_3$};
\node[mycirc,label=above:{\tiny $[1,3]$},fill=black,text=white] (v2) at ( 60:3) {$v_2$};
\node[mycirc,label=above:{\tiny $[1,1]$},fill=black,text=white] (v1) at (2*60:3) {$v_1$};
\node[mycirc,label=above left:{\tiny $[2,2]$}] (v6) at (3*60:3) {$v_6$};
\node[mycirc,label=below:{\tiny $[3,5]$}] (v5) at (4*60:3) {$v_5$};
\node[mycirc,label=below:{\tiny $[5,5]$}] (v4) at (5*60:3) {$v_4$};

\node[mytext] (e) at (270:6) {(e)};

    \draw[<->] (v1) -- (v2);
    \draw[<->, line width=1pt] (v1) -- (v6);
    \draw[<->, line width=1pt] (v6) -- (v5);
    \draw[<->, line width=1pt] (v4) -- (v5);
    \draw[<->] (v3) -- (v2);
    \draw[<->] (v2) -- (v6);
    \draw[<->, line width=1pt] (v4) -- (v3);
    
     \draw[<->,dashed] (v2) -- (v5);
       \draw[<->,dashed] (v3) -- (v5);
\end{scope}

\end{tikzpicture}
\caption{An example of a network and its critical additive edge-sets.}
\label{critadd}
\end{figure}

In Fig. \ref{critadd} (a), a graph $G$ with a black  ZFS is shown. In Fig. \ref{critadd} (b)-(e), different sets of chains $\mathcal{C}$ and the time intervals $\mathbb{T}(v)=[T(v), T_{\mathrm{max}}(v)]$ for every $v\in V(G)$ are given, and the associated perfect graphs are shown (only the bidirectional edges of the perfect graphs are shown for clarity). By Theorem \ref{th5}, $n^c_{\mathrm{add}}(G)=16$. 
In Fig. \ref{critadd} (b),
 a critical  edge-set is $E^c_{\mathrm{add}}(G)=\{(v_3,v_6), (v_6,v_3), (v_4,v_6),(v_6,v_4),(v_1,v_1),(v_2,v_2),(v_3,v_3),(v_4,\\v_4),(v_5,v_5),(v_6,v_6),(v_3,v_1),(v_4,v_1),(v_5,v_1),(v_4,v_2),(v_5,v_2),\\(v_5,v_3)\}$.


Next, we study the robustness of ss-controllability of a network with respect to edge removals and describe a critical subtractive edge-set. A formula for the critical 
subtractive number is then presented.

\color{black}\begin{pro}
\color{black}Consider a network with the graph $G$ and the ZFS $V_C$. If $G\in\mathcal{G}^{\mathcal{C},T}$, for a set of node-disjoint chains $\mathcal{C}$ with the set of sources $V_C$ and  time function $T$, then $E(G)\setminus \bigcup_{i=1}^m E(C_i) $ is a critical subtractive edge-set of $G$. Moreover, $n^c_{\mathrm{sub}}(G)=|E(G)|-n+m$. 
\label{th7}
\end{pro}

\color{black}\textit{Proof:} By Definition \ref{de3}, $E'=E(G)\setminus \bigcup_{i=1}^m E(C_i) $ includes edges which may or may not exist in a $(\mathcal{C},T)$-constructed graph. Then, if we remove any subset of $E'$ from the edge set of $G$, we obtain a new $(\mathcal{C},T)$-constructed graph which is still SSC  from Theorems 1 and \ref{th2}. Moreover, $E'$ is the largest set of edges which can be removed from the edge set of a network so that its ss-controllability is preserved. Because, by removing more than $|E'|$ edges from the network, a graph with more than $m$ connected components is obtained, which cannot be SSC from only $m$ control nodes.  Hence, $E'$ is a critical subtractive edge-set of $G$. Moreover, since $\sum_{i=1}^m |E(C_i)|=n-m$, we have $n^c_{\mathrm{sub}}(G)=|E(G)|-n+m$.
\carre

\section{Network combinations}

In this section, methods of combination of networks, resulting in \color{black}a SSC \color{black} network-of-networks are presented.

\subsection{Combining SSC networks} 
We now present methods for combining SSC networks while preserving their ss-controllability. We also consider a structural uncertainty for the networks and propose 
methods for combining them such that despite the uncertainty in their respective structures, 
the combined network remains SSC. 

For any $1\leq i\leq l$, consider an LTI network with graph $G_i=(V(G_i),E(G_i))$ of size $n_i$, that is SSC. Let $V_C^i$ be the set of control nodes, and $|V_C^i|=m_i$. Consider a set of node-disjoint chains $\mathcal{C}^i=\{C^i_1,\ldots,C^i_{m_i}\}$ and a time function $T^i$ such that $G_i\in\mathcal{G}^{\mathcal{C}^i,T^i}$.  

Now, consider a graph $G=(V(G),E(G))$, defined as $G=\mathrm{comb}(\bigcup_{i=1}^lG_i,E^*)$, which implies that $G$ is a combination of the graphs $G_1,\ldots, G_l$     through connecting some of their nodes via directed edges in $E^*$. Then, we have $V(G)=\bigcup_{i=1}^l V(G_i)$, and $E(G)=\bigcup_{i=1}^l E(G_i)\cup E^*$, where for any $(u,v)\in E^*$, there is some $1\leq k,j\leq l$, such that $k\neq j$, and $u\in V(G_k)$, and $v\in V(G_j)$. Thus, $E^*$ includes all edges that are between nodes of   two different graphs. 

Let us  define $\mathcal{G}_u$ as the class of all graphs $ G'=\mathrm{comb}(\bigcup_{i=1}^lG'_i,E^*)$, where $G'_i\in\mathcal{G}^{\mathcal{C}^i,T^i}$, for $i=1,\ldots,l$.
In fact,  one can consider structural uncertainties in any $G_i$ and obtain a graph in $\mathcal{G}^{\mathcal{C}^i,T^i}$; then, by combining these graphs via edges in $E^*$, a graph in $\mathcal{G}_u$ is provided. It is obvious that $G\in \mathcal{G}_u$.

 We aim to propose a method to combine graphs $G_1,\ldots,G_l$ in a way that: (1) the LTI network with graph $G=\mathrm{comb}(\bigcup_{i=1}^lG_i,E^*)$ is SSC, and (2) \color{black} all the LTI networks with graphs in $\mathcal{G}_u$  are SSC. \color{black}

Given a time function $T^i$ for every graph $G_i$, $i=1,\ldots,l$, Algorithm 2 transcribes an update on $T^i$ for the purpose of combining networks. It is noted that in Algorithm 2, for a sequence $S$, $|S|$ is the number of its elements, and $S(j)$ denotes its $j$th element. 

\textit{Procedure 1:} Consider  graphs $G_1,\ldots,G_l$, where $G_i\in\mathcal{G}^{\mathcal{C}^i,T^i}$, $i=1,\ldots,l$. Also, let $S$ be a sequence of $G_i$'s, $i=1,\ldots,l$, where every $G_i$ is repeated $q_i=n_i-m_i$ times.  The number of different sequences which can be made in this way is $(\sum_{i=1}^l q_i)!/ (q_1!\ldots q_l!)$. Considering the sequence $S$, run Algorithm 2 for any of the  graphs $G_1,\ldots, G_l$ ($G^*$ in Algorithm 2 can be any of the graphs $G_1,\ldots,G_l$), and obtain a new time function $T:V(G)\rightarrow [1,n-m+1]$ for $G$, where $n=\sum_{i=1}^l n_i$, and $m=\sum_{i=1}^l m_i$. Moreover, let $\mathcal{C}=\{\mathcal{C}^1,\ldots,\mathcal{C}^l\}$.  Note that the  computational complexity of Algorithm 2 is $\mathcal{O}(n)$. 

\begin{figure}[htb]
\mbox{}\hrulefill
\vspace{-.4em}
\\
\textbf{Algorithm 2:}

\vspace{-.7em}
\mbox{}\hrulefill
\vspace{-.01em}
\\
\small Given a sequence $S$ and the graph $G^*$ with the time function $T^*$; \\
j=0;\\
For $k=1:|S|$\\
\begin{tabular}{|l}
If $S(k)==G^*$\\
\begin{tabular}{|l}
$j=j+1$;\\
For $v\in V(G^*)$ that $T^*(v)=j+1$;\\
\begin{tabular}{|l}
$T(v)=k+1$;
\end{tabular}\\
end for\\
\end{tabular}\\
end if\\
\end{tabular}\\
end for\\
\vspace{-.5em}
\mbox{}\hrulefill
\caption{An algorithm that with a given sequence $S$  and a time function $T^*$, updates the integer assigned to every node of a graph $G^*$.}
\vspace{-0.08in}
\end{figure}
For example, consider  the graphs $G_1$ and $G_2$ in Fig. \ref{combin} (a). The time interval $\mathbb{T}^1(v)=[T^1(v), T^1_{\mathrm{max}}(v)]$ associated with every node $v\in V(G_1)$ and the time interval $\mathbb{T}^2(v)=[T^2(v), T^2_{\mathrm{max}}(v)]$ associated with every node $v\in V(G_2)$ are also given. Now, consider a sequence $S$ of $G_1$ and $G_2$, where $G_i$ ($i=1,2$) is repeated $n_i-m_i=2$ times. For example, let $S=(G_2,G_1,G_1,G_2)$. In Fig. \ref{combin} (b), the updated time intervals assigned to the  nodes of the combined network, obtained by running Algorithm 2, are presented. \color{black}  \color{black}

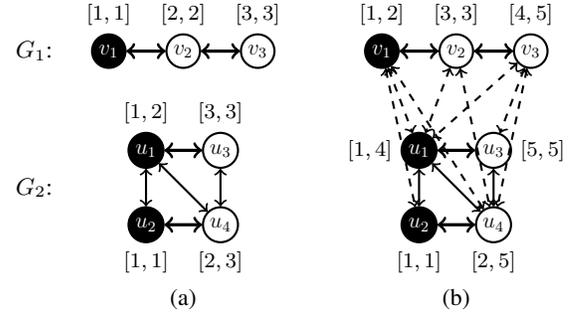
\begin{figure}[!h]
\centering
\begin{tikzpicture}[scale=.33]
\tikzset{mycirc/.style={draw,shape=circle,style=thick,inner sep=2pt,scale=.7}}
\tikzset{mytext/.style={shape=circle,style=thick}}
    ]
    \node[mycirc,label=above:{\small $[1,1]$},fill=black,text=white] (v1) at (0,0) {\Large $v_1$};
    \node[mycirc,label=above:{\small $[2,2]$}] (v2) at (3,0) {\Large $v_2$};
    \node[mycirc,label=above:{\small $[3,3]$}] (v3) at (6,0) {\Large $v_3$};
\node[mytext] (G1) at (-3,0) {$G_1$:};  
\node[mytext] (a) at (3,-10) {(a)};

    \draw[<->, style=very thick] (v1) -- (v2);
    \draw[<->, style=very thick] (v2) -- (v3);

    \node[mycirc,label=above:{\small $[1,2]$},fill=black,text=white] (u1) at (1.5,-4) {\Large $u_1$};
    \node[mycirc,label=below:{\small $[1,1]$},fill=black,text=white] (u2) at (1.5,-7) {\Large $u_2$};
    \node[mycirc,label=above:{\small $[3,3]$}] (u3) at (4.5,-4) {\Large $u_3$};
    \node[mycirc,label=below:{\small $[2,3]$}] (u4) at (4.5,-7) {\Large $u_4$};
    
    \node[mytext] (G2) at (-3,-5.5) {$G_2$:}; 
    
    \draw[<->, style=thick] (u1) -- (u2);
    \draw[<->, style=very thick] (u2) -- (u4);
    \draw[<->, style=very thick] (u1) -- (u3);
    \draw[<->, style=thick] (u4) -- (u3);
    \draw[<->, style=thick] (u1) -- (u4);



\begin{scope}[shift={(11,0)}]

\tikzset{mycirc/.style={draw,shape=circle,style=thick,inner sep=2pt,scale=.7}}
\tikzset{mytext/.style={shape=circle,style=thick}}
    ]
    \node[mycirc,label=above:{\small $[1,2]$},fill=black,text=white] (v1) at (0,0) {\Large $v_1$};
    \node[mycirc,label=above:{\small $[3,3]$}] (v2) at (3,0) {\Large $v_2$};
    \node[mycirc,label=above:{\small $[4,5]$}] (v3) at (6,0) {\Large $v_3$};
    \node[mytext] (b) at (3,-10) {(b)};    

    \draw[<->, style=very thick] (v1) -- (v2);
    \draw[<->, style=very thick] (v2) -- (v3);

    \node[mycirc,label=left:{\small $[1,4]$},fill=black,text=white] (u1) at (1.5,-4) {\Large $u_1$};
    \node[mycirc,label=below:{\small $[1,1]$},fill=black,text=white] (u2) at (1.5,-7) {\Large $u_2$};
    \node[mycirc,label=right:{\small $[5,5]$}] (u3) at (4.5,-4) {\Large $u_3$};
    \node[mycirc,label=below:{\small $[2,5]$}] (u4) at (4.5,-7) {\Large $u_4$};
    
    
    \draw[<->, style=thick] (u1) -- (u2);
    \draw[<->, style=very thick] (u2) -- (u4);
    \draw[<->, style=very thick] (u1) -- (u3);
    \draw[<->, style= thick] (u4) -- (u3);
    \draw[<->, style= thick] (u1) -- (u4);

    \draw[<->, dashed, thick] (v1) -- (u2);
    \draw[<->, dashed, thick] (v1) -- (u1);
    \draw[<->, dashed,thick] (v1) -- (u4);
    \draw[<->, dashed,thick] (v2) -- (u1);
    \draw[<->, dashed,thick] (v2) -- (u4);
    \draw[<->, dashed,thick] (v3) -- (u1);
    \draw[<->, dashed,thick] (v3) -- (u3);
    \draw[<->, dashed,thick] (v3) -- (u4);
\end{scope}
\end{tikzpicture}
\caption{a) Graphs $G_1$ and $G_2$, b) their combination.}
\label{combin}
\end{figure}

\begin{theo}
Consider an LTI network with a graph $G=\mathrm{comb}(\bigcup_{i=1}^lG_i,E^*)$, where $G_i\in\mathcal{G}^{\mathcal{C}^i,T^i}$, for $i=1,\ldots,l$.  
 Let a set of node-disjoint chains $\mathcal{C}$ and a time function  $T$ be provided by Procedure 1. Then, any network with a graph $G\in\mathcal{G}_u$  is  SSC if the following condition holds: for all $u\in V(G_i) $ and $v\in V(G_j)$ ($i\neq j$), $(u,v)\notin E^*$ if $T_{\mathrm{max}}(u)<T(v)$. Moreover, the largest set of edges which can be added between graphs while preserving the ss-controllability of any network with a graph in $\mathcal{G}_u$ is $E^*_{\max}=\{(u,v):u\in V(G_i), v\in V(G_j), 1\leq i,j\leq l, \:\:i\neq j, \:\:\mbox{and}\:\:T_{\mathrm{max}}(u)\geq T(v) \}$, and  $|E^*_{\max}|= \frac{1}{2} n( n+1)+\frac{1}{2} m(2n-m-1) -\sum_{i=1}^l\left( \frac{1}{2} n_i(n_i+1)+\frac{1}{2} m_i(2n_i-m_i-1) \right)$.
\label{th8}     
\end{theo}

\textit{Proof:} To show the ss-controllability, 
it suffices to prove that $G$ is a $(\mathcal{C}, T)$-constructed graph. By the assumption, if $i\neq j$ and $u\in V(G_i) $ and $v\in V(G_j)$, we have $(u,v)\notin E(G)$ if $T_{\mathrm{max}}(u)<T(v)$. Then, we should only prove that for $i=1,\ldots, l$ and for  all $u,v \in V(G_i)$ that $(u,v)\notin \bigcup_{k=1}^{m^i}E(C_k^i)$, one has $(u,v)\notin E(G)$ if $T_{\mathrm{max}}(u)<T(v)$. Note that since $G_i\in \mathcal{G}^{\mathcal{C}^i, T^i}$, by Definition \ref{de3}$, (u,v)\notin E(G_i)$ if $T^i_{\mathrm{max}}(u)<T^i(v)$. Now, we claim that $T^i_{\mathrm{max}}(u)<T^i(v)$ if and only if $T_{\mathrm{max}}(u)<T(v)$. Since $T^i_{\mathrm{max}}(u)=T^i(u+1)-1$ and $T_{\mathrm{max}}(u)=T(u+1)-1$, we should prove that for every $w,v\in V(G_i)$, $i=1,\ldots,m$, $T^i(w)\leq T^i(v)$ if and only if $T(w)\leq T(v)$. From Algorithm 2, one can see that if for some $u\in V(G_i)$, $T^i(u)=j+1$, and $T(u)=k+1$, then $j$th $G_i$ in the sequence $S$ is the $k$th element of $S$. 
 Now, let $j_1$th and $j_2$th $G_k$ in the sequence $S$ be respectively, its $k_1$th and $k_2$th elements. Then, one can see that $j_1\leq j_2$ 
if and only if $k_1\leq k_2$. Accordingly, $T^i(w)\leq T^i(v)$ if and only if $T(w)\leq T(v)$, and $G$ is a $(\mathcal{C}, T)$-constructed graph. 
From Theorem \ref{th3}, a SSC network with graph $G$ has  maximum number of edges when $G$ is a perfect $(\mathcal{C},T)$-constructed graph. When $G_i$'s, $i=1,\ldots,l$ are all perfect, the maximal set of edges which can be added  is $E^*_{\max}$ with a cardinality  obtained from Lemma \ref{th4}.    
\carre

\subsection{Combination of directed acyclic networks}
We now propose a method for combining networks with acyclic structures such that the 
corresponding network-of-networks is SSC with only a single control node. 
A directed acyclic graph is a directed graph with no directed cycles. However, if all edges are replaced with undirected ones, there may be some undirected cycles in the graph. 

\begin{deff}
A topological ordering of a directed graph $G$ is a linear ordering of the nodes such that for every $(u,v)\in E(G)$, $v$ comes before $u$ in the ordering. Hence, where there is a topological ordering for a graph $G$, one can index the nodes in a way such that for all $(i,j)\in E(G)$, $i>j$.
\label{de7}
\end{deff}

\begin{lem}[\cite{bang2008digraphs}]
A graph has a topological ordering if and only if it is a directed acyclic graph. 
\label{le2}
\end{lem}
  
A directed acyclic graph has at least one topological ordering, but a topological ordering may not be unique. There are some algorithms computing a topological ordering of a graph. For example, in Kahn's algorithm it can be computed in $\mathcal{O}(|V|+|E|)$ \cite{kahn1962topological}.

Consider $l$  directed acyclic graphs $G_1,\ldots,G_l$. Let $n_i=|V(G_i)|$, and $n=\sum_{i=1}^l n_i$. Assume that the nodes of each graph are indexed by the topological ordering. The following graph composition procedure can lead to \color{black}a SSC \color{black} network with a single control node.

\textit{Procedure 2:} Arrange $G_1,\ldots,G_l$ in a sequence $S$ in a way that for all $i=1,\ldots,l$, $G_i$ is repeated $n_i$ times. Moreover, no two $G_i$'s are put beside each other in the sequence. Given the sequence $S$, for all $i=1,\ldots,l$, index $G_i$'s in $S$ by a superscript, considering their place in the sequence. \color{black}More precisely\color{black}, index the $j$th $G_i$ in $S$ as $G^j_i$, for $j=1,\ldots, n_i$. For example, assume that $l=3$, $n_1=2$, $n_2=3$, and $n_3=1$. Let $S=(G_1, G_2,G_3,G_2,G_1,G_2)$. Then, we can write $S=(G_1^1, G_2^1,G_3^1,G_2^2,G_1^2, G_2^3)$. Let $v^i$ be  node $v$ of the graph $G_i$. Now, for all $k=1,\ldots, n-1$, if $S(k)=G_i^{v}$ and $S(k+1)=G_j^{u}$, for some $1\leq i, j\leq l$, $1\leq v\leq n_i$, and $1\leq u\leq n_j$, then add an edge from the node $v^i$ to the node $u^j$. Let $G$ be the obtained graph, that is, a combination of graphs $G_1,\ldots,G_l$. Now, 
if  $S(k)=G_i^{v}$, let $T(v^i)=k$, where $T$ is a time function defined on nodes of $G$.


\begin{pro}
Let $v$ be a node that $T(v)=1$ and $V_C=\{v\}$. An LTI network with graph $G$, obtained by Procedure 2, is SSC. 
\label{pr3}
\end{pro}    
   
\textit{Proof:} Let $C_1$ be a chain with $V(C_1)=\bigcup_{i=1}^l V(G_i)$, and $E(C_1)=\{(u,v): T(u)=k, T(v)=k+1, k=1\ldots, n-1\}$. Define $\mathcal{C}=\{C\}$. We claim $G$ is a $(\mathcal{C},T)$-constructed graph. Otherwise, for some $u, v\in V(G)$ that $T_{\mathrm{max}}(u)<T(v)$ and $(u,v)\notin E(C)$, we have $(u,v)\in E(G)$. Then, there is  some graph $G_i$, where $u, v\in V(G_i)$. Moreover, since $T_{\mathrm{max}}(u)<T(v)$, then $u<v$. Thus, by Lemma \ref{le2} and Definition \ref{de7}, there cannot be an edge of the form $(u,v)$ in $G_i$, which is a contradiction. Then, $G\in \mathcal{G}^{\mathcal{C},T}$, and the source of $\mathcal{C}$ renders it SSC. 
\carre

\color{black}
\section{Controllability of a family of LTV networks}   

As shown in previous sections, ss-controllability of an LTI network can remain intact while there is an uncertainty in the presence or absence of certain edges. In this section, we present conditions for ss-controllablity of an LTV network, where not only edge weights can be time-varying, but also over some time intervals, certain edges can be removed from, or added to, the network. 

For a given graph $G$, consider a set of time-varying $n\times n$ matrices $A(t)$, where for \emph{almost every} (a.e.) $t\in \mathbb{R}$, $A(t)\in \mathcal{Q}(G)$. Thus, for {a.e.} $t\in \mathbb{R}$ and $i\neq j$, we have $A_{ij}(t)\neq 0$ if and only if $(j,i)\in E(G)$.
Note that given an edge $(j,i)\in E(G)$, there might exist some $t$ for which $A_{ij}(t)=0$, but the  (Lebesgue) measure of such a set of $t$'s is zero. That is, there is no time interval over which $A_{ij}(t)=0$. In \cite{reissig2014strong}, the following controllability result for such LTV networks has been presented. In subsequent discussion, let $t_0,t_1\in \mathbb{R}$, and $t_0<t_1$. 


 
 \begin{theo}
 An LTV network of the form (\ref{lv}), where for a.e. $t\in \mathbb{R}$, $A(t)\in \mathcal{Q}(G)$,  is controllable on $[t_0,t_1]_{\mathbb{R}}$  if and only if $V_C$ is a ZFS of $G$. 
\label{ltv1}
\end{theo}

Now,  given a set of node-disjoint chains $\mathcal{C}=\{C_1,\ldots,C_m\}$ and a time function $T$, consider the class of $(\mathcal{C},T)$-constructed graphs $\mathcal{G}^{\mathcal{C},T}$. Let $A(.)$ be a piecewise continuous function of $t$, and for {a.e.} $t\in \mathbb{R}$, $A(t)\in\mathcal{P}(\mathcal{G}^{\mathcal{C},T})$, that is, $A(t)\in \mathcal{Q}(G)$, for some $G\in \mathcal{G}^{\mathcal{C},T}$. For example, for a time interval $[t_0,t_1]_{\mathbb{R}}$, we might have $A(t)\in\mathcal{Q}(G_1) $, and for other time interval $[t_1,t_2]_{\mathbb{R}}$, $A(t)\in\mathcal{Q}(G_2) $, where $G_1, G_2\in \mathcal{G}^{\mathcal{C},T}$. If $(j,i)\in E(C_i)$, for some $1\leq i\leq m$, then we have $A_{ij}(t)\neq 0$, for a.e. $t\in \mathbb{R}$. Otherwise, we may have $(j,i)\in E(G_1)$, while $(j,i)\notin E(G_2)$, and then for a.e. $t\in[t_0,t_1]_{\mathbb{R}}$, $A_{ij}(t)\neq0$, while for a.e.  $t\in[t_1,t_2]_{\mathbb{R}}$, $A_{ij}(t)=0$. 

For example, consider a chain $C$ with $V(C)=\{v_1,v_2,v_3\}$, and $E(C)=\{(v_1,v_2),(v_2,v_3)\}$. Define $T(v_i)=i$, $1\leq i\leq 3$, and $\mathcal{C}=\{C\}$. Now, consider some $A(.)$ that for a.e. $t\in \mathbb{R}$, $A(t)\in \mathcal{P}(\mathcal{G}^{\mathcal{C},T})$.  Let $t_0=0$, and $t_1=10$.  Thus, for a.e. $t\in[0,10]_{\mathbb{R}}$, we must have $A_{21}(t)\neq 0$, $A_{32}(t)\neq 0$, and $A_{31}(t)=0$. Moreover, we may assume that $A_{11}(t)\neq 0$, for $t\in[0,2]_{\mathbb{R}}\cup [4,6]_{\mathbb{R}} $; $A_{12}(t)\neq 0$, for $t\in [1,3]_{\mathbb{R}}$; $A_{13}(t)\neq 0$, for $t\in [2,5]_{\mathbb{R}}$; $A_{22}(t)\neq 0$, for $t\in[7,10]_{\mathbb{R}}$; $A_{23}(t)\neq 0$, for $t\in [0,3]_{\mathbb{R}}\cup [5,6]_{\mathbb{R}}$; and $A_{33}(t)\neq 0$, for $t\in [4,7]_{\mathbb{R}}$.  

In the following, we establish a controllability condition for all networks with these system matrices.


  
  \begin{theo}
Given a class $\mathcal{G}^{\mathcal{C},T}$, any LTV system (\ref{lv}), where  for  a.e. $t\in \mathbb{R}$, $A(t)\in\mathcal{P}(\mathcal{G}^{\mathcal{C},T})$, is controllable on $[t_0,t_1]_{\mathbb{R}}$  if and only if $V_C$ is any set that includes the sources of $\mathcal{C}$. 
\label{ltv2}
\end{theo}
 
\emph{Proof:}
Let $\mathcal{S}$ be the set of sources of $\mathcal{C}$. First, consider systems whose system matrix $A(t)$ is  constant in time. Then, from Theorem \ref{p},  for the controllability of an LTI system with $A\in\mathcal{P}(\mathcal{G}^{\mathcal{C},T})$, we should have $\mathcal{S}\subseteq V_C$, and the necessary condition is proved.

Now, assume that $V_C=\mathcal{S}$. To prove the sufficiency, consider a system in this family with system matrix $A(t)$ and the transition matrix $\Phi(.,.)$. Since for a.e. $t\in \mathbb{R}$, $A(t)\in\mathcal{P}(\mathcal{G}^{\mathcal{C},T})$, then for any $(j,i)\in \bigcup_{i=1}^m E(C_i)$, we have $A_{ji}(t)\neq 0$, for a.e. $t\in[t_0,t_1]_{\mathbb{R}}$.

Let $t^0=t_0$ and $t^p=t_1$. Then, for some $p\geq 1$, one can find a time sequence $(t^0, t^1,\ldots,t^p)$, where $t^i<t^{i+1}$, such that the following property holds for all $i,j\in\{1,\ldots, n\}$: either  $A_{ij}(t)=0$ or $A_{ij}(t)\neq 0$, for a.e. $t\in[t^i,t^{i+1}]$, $1\leq i\leq p$. Hence for a.e. $t$ in any time sub-interval $[t^i,t^{i+1}]$, $A(t)$ is of the same zero-nonzero pattern. Then, there is a graph $G^i\in \mathcal{G}^{\mathcal{C},T}$, such that for a.e. $t\in [t^i,t^{i+1}]$, we have $A(t)\in \mathcal{Q}(G^i)$, $0\leq i\leq p-1$.  Now, consider an LTV system  with the transition matrix $\Phi^*(.,.)$ and the system matrix $A^*(t)$, where for a.e. $t\in\mathbb{R}$, $A^*(t)\in \mathcal{Q}(G^i)$, and   assume that for a.e. $t\in [t^i,t^{i+1}]$, we have $A^*(t)=A(t)$. 
From (\ref{tm}), one can see that  the transition matrix $\Phi(t_2,t_1)$ is a function that depends only on the numerical values of the matrix $A(t)$, for $t\in[t_1,t_2]$. 
 Thus, for all $\tau\in [t^i,t^{i+1}]$, we have $\Phi(t^{i+1},\tau)=\Phi^*(t^{i+1},\tau)$. From Theorem \ref{ltv1}, one can see that since $G^i$ is a $(\mathcal{C},T)$-constructed graph, and  $\mathcal{S}$ is a ZFS of $G^i$, then the LTV network with the system matrix $A^*(t)$ is controllable on $[t^i,t^{i+1}]$. Thus, based on Proposition \ref{tm2}, the condition $\nu^T\Phi^*(t^{i+1},\tau)B=0$, for a.e. $\tau\in [t^i,t^{i+1}]$, implies that $\nu=0$. Accordingly, since for all $\tau\in [t^i,t^{i+1}]$, we have $\Phi^*(t^{i+1},\tau)=\Phi(t^{i+1},\tau)$, then $\nu^T\Phi(t^{i+1},\tau)B=0$ implies $\nu=0$. Thus, from Proposition \ref{tm2}, for all $1\leq i\leq p-1$, the network with system matrix $A(t)$ is controllable on $[t^i,t^{i+1}]$. That is, for any $1\leq i\leq p-1$, the state of the system can be driven from any initial state at time $t^i$ to any final state at time $t^{i+1}$; similarly, from any initial state at $t_0$ to any final state at $t_1$, completing the proof.
 %
\carre

Hence, we have extended  Theorem \ref{ltv1} to the LTV networks, that in addition to their edge weights, 
have a time-varying structure. 

\color{black} 
\section{Conclusion}

In this paper, we examined the preservation of ss-controllability of an LTI network under structural perturbations. We then proceeded to show that if the number of edges added to, or removed from a network is bigger than certain bounds, ss-controllability is destroyed. These bounds depend only on the size of the  network and the number of its control nodes. Moreover, we described the maximal sets of edges, adding or removing any subset of which, do not disturb the network ss-controllability. Furthermore, we provided combinatorial algorithms for combining networks in order to ensure \color{black}a SSC \color{black} network-of-networks.  Finally, we derived controllability conditions for a family of LTV networks with time-varying structures.       

\bibliographystyle{IEEEtran}
\bibliography{library}

\begin{thebibliography}{10}
\providecommand{\url}[1]{#1}
\csname url@samestyle\endcsname
\providecommand{\newblock}{\relax}
\providecommand{\bibinfo}[2]{#2}
\providecommand{\BIBentrySTDinterwordspacing}{\spaceskip=0pt\relax}
\providecommand{\BIBentryALTinterwordstretchfactor}{4}
\providecommand{\BIBentryALTinterwordspacing}{\spaceskip=\fontdimen2\font plus
\BIBentryALTinterwordstretchfactor\fontdimen3\font minus
  \fontdimen4\font\relax}
\providecommand{\BIBforeignlanguage}[2]{{%
\expandafter\ifx\csname l@#1\endcsname\relax
\typeout{** WARNING: IEEEtran.bst: No hyphenation pattern has been}%
\typeout{** loaded for the language `#1'. Using the pattern for}%
\typeout{** the default language instead.}%
\else
\language=\csname l@#1\endcsname
\fi
#2}}
\providecommand{\BIBdecl}{\relax}
\BIBdecl

\bibitem{mousavi2017robust}
S.~S. Mousavi, M.~Haeri, and M.~Mesbahi, ``Robust strong structural
  controllability of networks with respect to edge additions and deletions,''
  in \emph{Proc. American Control Conf.}, Seattle, 2017, pp. 5007--5012.

\bibitem{mousavi2018controllability}
------, ``Controllability analysis of threshold graphs and cographs,'' in
  \emph{Proc. European Control Conf.}, Limassol, Cypres, 2018, pp. 1--6.

\bibitem{mousavi2018laplacian}
------, ``Laplacian dynamics on cographs: Controllability analysis through
  joins and unions,'' \emph{arXiv:1802.03599}, 2018.

\bibitem{lin1974structural}
C.-T. Lin, ``Structural controllability,'' \emph{IEEE Trans. Automat. Contr.},
  vol.~19, no.~3, pp. 201--208, 1974.

\bibitem{liu2011controllability}
Y.-Y. Liu, J.-J. Slotine, and A.-L. Barab{\'a}si, ``Controllability of complex
  networks,'' \emph{Nature}, vol. 473, no. 7346, pp. 167--173, 2011.

\bibitem{mayeda1979strong}
H.~Mayeda and T.~Yamada, ``Strong structural controllability,'' \emph{SIAM J.
  Contr. and Optimiz.}, vol.~17, no.~1, pp. 123--138, 1979.

\bibitem{bowden2012strong}
C.~Bowden, W.~Holderbaum, and V.~M. Becerra, ``Strong structural
  controllability and the multilink inverted pendulum,'' \emph{IEEE Trans.
  Automat. Contr.}, vol.~57, no.~11, pp. 2891--2896, 2012.

\bibitem{jarczyk2011strong}
J.~C. Jarczyk, F.~Svaricek, and B.~Alt, ``Strong structural controllability of
  linear systems revisited,'' in \emph{Proc. 50th {IEEE} Conf. on Decision and
  Control and Eur. Control Conf.}, Orlando, FL, 2011, pp. 1213--1218.

\bibitem{chapman2013strong}
A.~Chapman and M.~Mesbahi, ``On strong structural controllability of networked
  systems: a constrained matching approach,'' in \emph{Proc. American Control
  Conf.}, Washington, DC, 2013, pp. 6126--6131.

\bibitem{monshizadeh2014zero}
N.~Monshizadeh, S.~Zhang, and M.~K. Camlibel, ``Zero forcing sets and
  controllability of dynamical systems defined on graphs,'' \emph{IEEE Trans.
  Automat. Contr.}, vol.~59, no.~9, pp. 2562--2567, 2014.

\bibitem{trefois2015zero}
M.~Trefois and J.-C. Delvenne, ``Zero forcing number, constrained matchings and
  strong structural controllability,'' \emph{Linear Alg. and its Applic.}, vol.
  484, pp. 199--218, 2015.

\bibitem{monshizadeh2015strong}
N.~Monshizadeh, M.~K. Camlibel, and H.~L. Trentelman, ``Strong targeted
  controllability of dynamical networks,'' in \emph{Proc. 54th {IEEE} Conf. on
  Decision and Control}, Osaka, 2015, pp. 4782--4787.

\bibitem{mousavi2016controllability}
S.~S. Mousavi and M.~Haeri, ``Controllability analysis of networks through
  their topologies,'' in \emph{Proc. 55th {IEEE} Conf. on Decision and
  Control}, Las Vegas, 2016, pp. 4346--4351.

\bibitem{van2017distance}
H.~J. van Waarde, M.~K. Camlibel, and H.~L. Trentelman, ``A distance-based
  approach to strong target control of dynamical networks,'' \emph{IEEE Trans.
  Automat. Contr.}, vol.~62, no.~12, pp. 6266--6277, 2017.

\bibitem{mousavi2018structural}
S.~S. Mousavi, M.~Haeri, and M.~Mesbahi, ``On the structural and strong
  structural controllability of undirected networks,'' \emph{IEEE Trans.
  Automat. Contr.}, vol.~63, no.~7, pp. 2234--2241, 2018.

\bibitem{mousavi2018null}
S.~S. Mousavi, A.~Chapman, M.~Haeri, and M.~Mesbahi, ``Null space strong
  structural controllability via skew zero forcing sets,'' in \emph{Proc.
  European Control Conf.}, Limassol, Cypres, 2018, pp. 1845--1850.

\bibitem{reissig2014strong}
G.~Reissig, C.~Hartung, and F.~Svaricek, ``Strong structural controllability
  and observability of linear time-varying systems,'' \emph{IEEE Trans.
  Automat. Contr.}, vol.~59, no.~11, pp. 3087--3092, 2014.

\bibitem{gracy2017structural}
S.~Gracy, F.~Garin, and A.~Y. Kibangou, ``Structural and strongly structural
  input and state observability of linear network systems,'' \emph{IEEE Trans.
  Control Netw. Syst.}, vol.~5, no.~4, pp. 2062--2072, 2018.

\bibitem{jia2018sufficient}
J.~Jia, H.~Trentelman, W.~Baar, and M.~Camlibel, ``A sufficient condition for
  colored strong structural controllability of networks,''
  \emph{IFAC-PapersOnLine}, vol.~51, no.~23, pp. 16--21, 2018.

\bibitem{sadat2019strong}
S.~S. Mousavi, M.~Haeri, and M.~Mesbahi, ``Strong structural controllability of
  signed networks,'' \emph{arXiv:1908.05732}, 2019.

\bibitem{rahimian2013structural}
M.~A. Rahimian and A.~G. Aghdam, ``Structural controllability of multi-agent
  networks: Robustness against simultaneous failures,'' \emph{Automatica},
  vol.~49, no.~11, pp. 3149--3157, 2013.

\bibitem{nie2014robustness}
S.~Nie, X.~Wang, H.~Zhang, Q.~Li, and B.~Wang, ``Robustness of controllability
  for networks based on edge-attack,'' \emph{PloS One}, vol.~9, no.~2, p.
  e89066, 2014.

\bibitem{mengiste2015effect}
S.~A. Mengiste, A.~Aertsen, and A.~Kumar, ``Effect of edge pruning on
  structural controllability and observability of complex networks,''
  \emph{Scientific Rep.}, vol.~5, 2015.

\bibitem{lu2016attack}
Z.-M. Lu and X.-F. Li, ``Attack vulnerability of network controllability,''
  \emph{PloS One}, vol.~11, no.~9, p. e0162289, 2016.

\bibitem{monnot2016sensitivity}
B.~Monnot and J.~Ruths, ``Sensitivity of network controllability to
  weight-based edge thresholding,'' in \emph{Complex Networks VII}, 2016, pp.
  45--61.

\bibitem{work2008zero}
{AIM Minimum Rank--Special Graphs Work Group}, ``Zero forcing sets and the
  minimum rank of graphs,'' \emph{Linear Alg. and its Applic.}, vol. 428,
  no.~7, pp. 1628--1648, 2008.

\bibitem{barioli2009minimum}
F.~Barioli, S.~M. Fallat, H.~T. Hall, D.~Hershkowitz, L.~Hogben, H.~Van~der
  Holst, and B.~Shader, ``On the minimum rank of not necessarily symmetric
  matrices: a preliminary study,'' \emph{Electron. J. Linear Algebra}, vol.~18,
  no.~1, pp. 126--145, 2009.

\bibitem{edholm2012vertex}
C.~J. Edholm, L.~Hogben, J.~LaGrange, and D.~D. Row, ``Vertex and edge spread
  of zero forcing number, maximum nullity, and minimum rank of a graph,''
  \emph{Linear Alg. and its Applic.}, vol. 436, no.~12, pp. 4352--4372, 2012.

\bibitem{berliner2013minimum}
A.~Berliner, M.~Catral, L.~Hogben, K.~Lied, and M.~Young, ``Minimum rank,
  maximum nullity, and zero forcing number of simple digraphs,''
  \emph{Electronic J. Linear Alg.}, vol.~26, no.~1, p.~52, 2013.

\bibitem{chapman2014controllability}
A.~Chapman, M.~Nabi-Abdolyousefi, and M.~Mesbahi, ``Controllability and
  observability of network-of-networks via cartesian products,'' \emph{IEEE
  Trans. Automat. Contr.}, vol.~59, no.~10, pp. 2668--2679, 2014.

\bibitem{yazicioglu2013leader}
A.~Y. Yazicioglu and M.~Egerstedt, ``Leader selection and network assembly for
  controllability of leader-follower networks,'' in \emph{Proc. American
  Control Conf.}, Washington, DC, 2013, pp. 3802--3807.

\bibitem{ji2016design}
Z.~Ji, T.~Chen, and H.~Yu, ``A design method for controllable topologies of
  multi-agent networks,'' in \emph{Proc. 35th Chinese Control Conf.}, Chengdu,
  China, 2016, pp. 7578--7584.

\bibitem{lu2015exploring}
C.~Lu, J.~Yu, R.-H. Li, and H.~Wei, ``Exploring hierarchies in online social
  networks,'' \emph{IEEE Trans. Knowledge and Data Engin.}, vol.~28, no.~8, pp.
  2086--2100, 2016.

\bibitem{elwert2013graphical}
F.~Elwert, ``Graphical causal models,'' in \emph{Handbook of Causal Analysis
  for Social Research}.\hskip 1em plus 0.5em minus 0.4em\relax Springer, 2013,
  pp. 245--273.

\bibitem{hespanha2018linear}
J.~P. Hespanha, \emph{Linear Systems Theory}.\hskip 1em plus 0.5em minus
  0.4em\relax Princeton Univ. Press, 2018.

\bibitem{sontag2013mathematical}
E.~D. Sontag, \emph{Mathematical Control Theory: Deterministic Finite
  Dimensional Systems}.\hskip 1em plus 0.5em minus 0.4em\relax New York:
  Springer Verlag, 1998.

\bibitem{hogben2012propagation}
L.~Hogben, N.~Kingsley, S.~Meyer, S.~Walker, and M.~Young, ``Propagation time
  for zero forcing on a graph,'' \emph{Disc. Applied Math.}, vol. 160, no.
  13-14, pp. 1994--2005, 2012.

\bibitem{kenter2018error}
F.~H. Kenter and J.~C.-H. Lin, ``On the error of a priori sampling: zero
  forcing sets and propagation time,'' \emph{Linear Alg. and its Applic.}, vol.
  576, pp. 124--141, 2019.

\bibitem{brimkov2019computational}
B.~Brimkov, C.~C. Fast, and I.~V. Hicks, ``Computational approaches for zero
  forcing and related problems,'' \emph{European J. Operational Research}, vol.
  273, no.~3, pp. 889--903, 2019.

\bibitem{bang2008digraphs}
J.~Bang-Jensen and G.~Z. Gutin, \emph{Digraphs: Theory, Algorithms and
  Applications}.\hskip 1em plus 0.5em minus 0.4em\relax Springer Science \&
  Business Media, 2008.

\bibitem{kahn1962topological}
A.~B. Kahn, ``Topological sorting of large networks,'' \emph{Communications of
  the ACM}, vol.~5, no.~11, pp. 558--562, 1962.

\end{thebibliography}

\end{document}